\documentclass[11pt,reqno]{amsproc}
\title{The vanishing surface tension limit of the Muskat problem}

\author{Patrick T. Flynn}
\address{Division of Applied Mathematics, Brown University, Providence, RI 02912}
\email{patrick\_flynn1@brown.edu}

\author{Huy Q. Nguyen}
\address{Department of Mathematics, Brown University, Providence, RI 02912}
\email{hnguyen@math.brown.edu}

\usepackage[margin=1in]{geometry}
\usepackage{amsmath, amsthm, amssymb, mathrsfs, stmaryrd}
\usepackage{times}
\usepackage{color}
\usepackage{hyperref}
\usepackage{graphicx}
\usepackage{subcaption}
\usepackage{float}

\newcommand{\bq}{\begin{equation}}
\newcommand{\eq}{\end{equation}}
\newcommand{\bqa}{\begin{eqnarray*}}
\newcommand{\eqa}{\end{eqnarray*}}

\theoremstyle{plain}
\newtheorem{theo}{Theorem}[section]
\newtheorem{prop}[theo]{Proposition}
\newtheorem{lemm}[theo]{Lemma}

\newtheorem{defi}[theo]{Definition}
\theoremstyle{definition}
\newtheorem{rema}[theo]{Remark}

\DeclareMathOperator{\dist}{dist}

\DeclareSymbolFont{pletters}{OT1}{cmr}{m}{sl}
\DeclareMathSymbol{s}{\mathalpha}{pletters}{`s}

\def\tt{\theta}
\def\eps{\varepsilon}
\def\na{\nabla}

\def\lb{\llbracket}
\def\rb{\rrbracket}
\def\les{\lesssim}

\def\mez{\frac{1}{2}}
\def\tdm{\frac{3}{2}}

\def\Rr{\mathbb{R}}
\def\T{\mathbb{T}}
\def\Nn{\mathbb{N}}

\def\cF{\mathcal{F}}

\def\cR{\mathcal{R}}
\def\ld{\lambda}

\def\p{\partial}
\def\na{\nabla}

\def\ka{\kappa}
\def\ma{\mathfrak{a}}

\def\wt{\widetilde}

\def\g{\mathfrak{g}}
\def\s{\mathfrak{s}}
\def\lbb{[\! [}
\def\rbb{]\! ]}

\numberwithin{equation}{section}

\date{today}
\begin{document}
\begin{abstract}
 The Muskat problem, in its general setting, concerns the interface evolution between two  incompressible fluids of different densities and viscosities in porous media.  The  interface motion is driven by gravity and capillarity forces, where the latter is due to surface tension. To leading order, both the Muskat problems with and without surface tension effect are  scaling invariant  in the Sobolev space $H^{1+\frac{d}{2}}(\Rr^d)$, where $d$ is the dimension of the interface. We prove that for  any subcritical data satisfying the Rayleigh-Taylor condition, solutions of the Muskat problem  with surface tension $\s$ converge to the unique solution of the Muskat problem without surface tension locally in time with the rate $\sqrt\s$  when $\s\to 0$. This allows for initial interfaces that have  unbounded or even not locally square integrable curvature. If in addition  the initial curvature is square integrable, we obtain the convergence with  optimal rate $\s$. 
\end{abstract}

\noindent

\maketitle

\section{Introduction}
In the studies of fluid flows, interfacial dynamics is broad and mathematically challenging. Some  interfacial problems in fluid dynamics that have been rigorously studied include the water wave problem, the compressible free-boundary Euler equations, the Hele-Shaw problem, the Muskat problem and the Stefan problem. The dynamics of the interface between the fluids strongly depends on  properties  of the fluids and of the media through which they flow. However,  a common feature in all the above problems is that the interface is driven by gravity and surface tension. Gravity is incorporated  in the momentum equations as an external force. On the other hand, surface tension balances the pressure jump across the interface (Young-Laplace equation):
\bq
\lb p\rb=\s H,
\eq
where $\lb p\rb$ is the pressure jump, $H$ is the mean-curvature of the interface, and $\s>0$ is the surface tension coefficient.  Well-posedness  in Sobolev spaces always holds when surface tension is taken into account but only holds under the Rayleigh-Taylor stability condition on the initial data when surface tension is neglected. 
It is a natural problem to justify the models without surface tension as the limit of the corresponding full models as surface tension vanishes. This question was addressed in \cite{AmbMas1, AmbMas2, ShaZen, CouHolShk, Ambrose2D, Ambrose3D, HadShk} for the problems listed above. The common theory is the following: if the initial data is stable, i.e. it satisfies  the Rayleigh-Taylor stability condition, and is {\it sufficiently smooth}, then solutions to the problem with surface tension converge to the unique solution of the problem without surface tension locally in time. The general strategy of proof consists of two points.
\begin{itemize}
\item[(i)] To leading order, the surface tension term $\s H$  provides a regularizing effect. For sufficiently smooth solutions, the difference between $\s H$ and its leading contribution can be controlled by the energy of the problem without surface tension. This yields a uniform time of existence $T_*$ for the problem with surface tension  $\s\to 0$.
\item[(ii)] For  sufficiently smooth solutions and for some $\tt\in [0, 1)$, the weighted mean curvature  $\s^\tt H$ is uniformly in $\s$ bounded (in some appropriate Sobolev norm) by the energy of the problem without surface tension. It follows that $\s H$, the difference between the two problems, vanishes as $\s^{1-\tt}$ as $\s\to 0$, establishing the convergence on the time interval $[0, T_*]$.  Note that the optimal rate corresponds to $\tt=0$.
\end{itemize}
Therefore, the vanishing surface tension limit becomes subtle if the initial data is sufficiently rough so that it can accommodate curvature singularities. As a matter of fact, in the aforementioned works, the initial curvature is at least bounded. In this paper, we prove that for the Muskat problem, the zero surface tension limit can be established for rough initial interfaces whose curvatures are {\it not bounded} or even {\it not locally $L^2$}. Regarding  quantitive properties of the zero surface tension limit,  the convergence rates in the aforementioned  works are either unspecified or suboptimal. In this paper, we  obtain the {\it optimal convergence rate} for the Muskat problem. The next subsections are devoted to a description of the Muskat problem and the statement of our main result. 

\subsection{The Muskat problem}
The Muskat problem \cite{Mus} concerns the interface evolution between two  fluids of densities $\rho^\pm$ and viscosities $\mu^\pm$ governed by Darcy's law for flows through porous media.  Specifically, the fluids occupy two domains $\Omega^\pm = \Omega^\pm_t \subset \mathbb R^{d+1}$ separated by an interface $\Sigma = \Sigma_t$, with $\Omega^+$ confined below a rigid boundary $\Gamma^+$, and $\Omega^-$ likewise above $\Gamma^-$.  We consider the case when the surfaces $\Gamma^\pm$ and $\Sigma$ are given by the graphs of functions, that is, we designate $b^\pm : \mathbb R^{d}_x \to \mathbb R$ and $\eta : \mathbb R_t\times \mathbb R^{d}_x \to \mathbb R$ for which
\begin{align}\label{domain}
&\Sigma = \{(x,\eta(t,x)) \ : \ x \in \mathbb R^d\},\\ 
&\Gamma^\pm= \{(x,b^\pm(x)) \ : \ x \in \mathbb R^{d} \},\\
&\Omega^- = \{(x,y) \in \mathbb R^{d}\times \mathbb R \ :  \ b^-(x) < y < \eta(t,x) \}, \\ 
&\Omega^+= \{(x,y) \in \mathbb R^{d}\times\mathbb R \ :   \eta(t,x) < y < b^+ (x) \},\\
&\Omega  = \Omega^+ \cup \Omega^-.
\end{align}
We also consider the case where one or both of $\Gamma^\pm = \emptyset$. In each domain $\Omega^\pm$, the fluid velocity $u^\pm$ and pressure $p^\pm$  obey Darcy's law:
\begin{equation}
\mu^\pm u^\pm + \nabla_{x,y} p^\pm  = -\rho^\pm g \vec{e}_{d+1},\quad \mathrm{div}_{x,y} u^\pm =0 \quad\text{ in } \Omega^{\pm},
\end{equation}
where $g$ denotes the gravitational acceleration, and $\vec{e}_{d+1}$ is the upward unit vector in the vertical direction.  For any two objects $A^+$ and $A^-$ associated with the domains $\Omega^+$ and $\Omega^-$ respectively, we denote the jump
\[
\lbb A \rbb = A^-- A^+
\]
whenever this difference is well-defined.  In particular, set
\[
\mathfrak g =g \lbb \rho \rbb.  
\]
 At the interface, there are three boundary conditions. First,  the
normal component of the fluid velocity is continuous across the interface
\begin{equation}
\lbb u \cdot n\rbb=0\quad\text{on}\quad \Sigma,
\end{equation}
where we fix $n$ to be the upward normal of the interface, specifically $n =\langle \nabla \eta\rangle^{-1} (-\nabla \eta, 1)$ with
\[
\langle \cdot \rangle = \sqrt{1+|\cdot|^2}.
\]
  Second, the interface is  transported by the normal fluid velocity, leading to the kinematic boundary condition
\begin{equation}
\eta_t = \langle \nabla \eta\rangle u^{-} \cdot n |_{\Sigma_t}.
\end{equation}
Third,  according to the Young-Laplace equation, the pressure jump is proportional to the mean curvature through surface tension:
\begin{equation}\label{bc:p}
\lbb p \rbb = \mathfrak s H(\eta) \equiv -\mathfrak s \text{div}( \langle \nabla \eta\rangle^{-1}\nabla \eta) \text{ on } \Sigma_t
\end{equation}
where $\mathfrak s\ge 0$ is the {\it surface tension coefficient} and 
\bq\label{def:H}
H(\eta)=-\text{div}( \langle \nabla \eta\rangle^{-1}\nabla \eta)
\eq
 is twice the mean curvature of $\Sigma$. 
 
Finally,  there is no transportation of fluid through the rigid boundaries:
\begin{equation}
u^\pm \cdot \nu ^\pm = 0\quad \text{on}\quad \Gamma^\pm,
\end{equation}
where $\nu^\pm = \pm\langle \nabla b\rangle^{-1}(-\nabla b^\pm, 1)$ is the outward normal of $\Gamma^\pm$.  If $\Gamma^\pm = \emptyset$, this condition is replaced by the decay condition
\bq
\lim_{y \to \pm \infty} u^\pm(x,y)=0.
\eq
 For the two-phase problem, we have $\rho^\pm$ and $\mu^\pm$ both as positive quantities. We will also consider the one-phase problem where the top fluid is treated as a vacuum by setting $\rho^+ = \mu^+ = 0$ and $\Gamma^+ = \emptyset$.

In the absence of the boundaries $\Gamma^\pm$, both the Muskat problems with and without surface tension  to leading order admit $\dot H^{1+\frac{d}{2}}(\Rr^d)$ as the scaling invariant Sobolev space in view of the scaling 
\[
\eta(x, t)\mapsto \ld^{-1}\eta(\ld x, \ld^3 t)\quad\text{and}\quad \eta(x, t)\mapsto \ld^{-1}\eta(\ld x, \ld t).
\]
In either case, the problem is quasilinear. The literature on well-posedness for the Muskat problem is vast. Early results can be found in  \cite{Chen, ConPug, EscSim, Yi, SieCafHow, Amb0,  Amb}. For more recent developments, we refer to \cite{CorGan07, CorCorGan, EscMat, CorCorGan, ConGanShvVic, GraShk}  for well-posedness, to \cite{ConCorGanStr, ConCorGanRodPStr, DenLeiLin, ConGanShvVic, GraShk, Cam, CorLaz, GanGarPatStr} for global existence, and to \cite{CasCorFefGan, CasCorFefGanMar} for singularity formation. Directly related to the problem addressed in the current paper is local well-posedness for low regularity large data. Consider first the problem without surface tension. In \cite{CheGraShk}, the authors obtained well-posedness for $H^2(\T)$ data for the one-phase problem, allowing for unbounded curvature. For the 2D Muskat problem without viscosity jump, i.e. $\mu^+=\mu^-$,  \cite{ConGanShvVic} proves well-posedness for data in all subcritical Sobolev spaces $W^{2, 1+}(\Rr)$. In $L^2$-based Sobolev spaces, \cite{Mat}  obtains well-posedness  for data in all subcritical spaces $H^{\tdm+}(\Rr)$. We also refer to \cite{AL} for a generalization of this result to homogeneous Sobolev spaces $\dot H^1(\Rr)\cap \dot H^s(\Rr)$, $s\in (\tdm, 2)$, allowing non-$L^2$ data. In \cite{NgPa}, local well-posedness for the Muskat problem in the general setting as described above was obtained for initial data in all subcritical Sobolev spaces $H^{1+\frac{d}{2}+}(\Rr^d)$, $d\ge 1$. The case of one fluid with infinite depth was independently obtained by \cite{AlaMeuSme}. Regarding the problem with surface tension, \cite{Mat, Mat2} consider initial data in $H^{2+}(\Rr)$. In the recent work \cite{Ng1}, well-posedness for data in all subcritical Sobolev spaces $ H^{1+\frac{d}{2}+}(\Rr^d)$ , $d\ge 1$, was established.
\subsection{Main Result}
In order to state the Rayleigh-Taylor stability condition solely in terms of the interface, we define the operator 
\bq
\mathrm{RT}(\eta)=1-\lbb  \mathfrak B(\eta) J (\eta)\rbb\eta\equiv 1-\Big(\mathfrak B^-(\eta) J^- (\eta)\eta-\mathfrak B^+(\eta) J^+ (\eta)\eta\Big),
\eq
where  $J^\pm(\eta)$ and  $\frak B^\pm(\eta)$ are respectively defined by and \eqref{def:J} and \eqref{Bop} below.
Our main result is  the following.
\begin{theo}\label{theo:main}
Consider either the one-phase Muskat problem or the two-phase Muskat problem in the stable regime $\rho^->\rho^+$.   The boundaries $\Gamma^\pm$ can be empty or graphs of  Lipschitz functions $b^\pm \in  \dot W^{1, \infty}(\Rr^d)$. For any $d\ge 1$, let $s > 1+\frac{d}{2}$ be an arbitrary  subcritical Sobolev index. Consider an initial datum $\eta_0\in H^s(\Rr^d)$ satisfying  
\begin{align}\label{cd:data:intro}
&\inf_{x\in \Rr^d}  \mathrm{RT}(\eta_0)\ge 2\ma>0, \quad \dist(\eta_0, \Gamma^\pm)\ge 2h>0.
\end{align}
Let $\s_n$ be a sequence of surface tension coefficients converging to $0$. Then, there exists $T_*>0$ depending only on $\| \eta_0\|_{H^s}$ and  $(\ma, h, s, \mu^\pm, \g)$ 
such that the following holds.

(i) The Muskat problems without surface tension and with surface tension $\s_n$ have a unique solution on $[0, T_*]$, denoted respectively by $\eta$ and $\eta_n$, that satisfy 
\begin{align}
&\eta_n\in C([0, T_*]; H^s(\Rr^d))\cap L^2([0, T_*]; H^{s+\tdm}(\Rr^d)),\\
&\eta\in L^\infty([0, T_*]; H^s(\Rr^d))\cap L^2([0, T_*]; H^{s+\mez}(\Rr^d))\cap C([0, T_*]; H^{s'}(\Rr^d))\quad\forall s'<s,\\ \label{uniformbound:theo}
 &\|(\eta_n, \eta)\|_{L^\infty([0,T_*];H^s)}+ \|(\eta_n, \eta)\|_{L^2([0,T_*];H^{s+\mez})} + \sqrt{\s} \|\eta_n\|_{L^2([0,T_*];H^{s+\tdm})}  \le  \mathcal F(\|\eta_0\|_{H^s}, \ma^{-1}),\\
 &\inf_{t\in [0, T_*]}\inf_{x\in \Rr^d}   \mathrm{RT}(\eta_n(t))>\ma, \quad \inf_{t\in [0, T_*]}\dist(\eta_n(t), \Gamma^\pm)>h,\\
 &\inf_{t\in [0, T_*]}\inf_{x\in \Rr^d}   \mathrm{RT}(\eta(t))>\ma, \quad \inf_{t\in [0, T_*]}\dist(\eta(t), \Gamma^\pm)> h,
 \end{align}
where $\cF:\Rr^+\times\Rr^+\to\Rr^+$  is nondecreasing and depends only on $(h, s, \mu^\pm, \g)$.

(ii)  As $n\to \infty$, $\eta_n$ converges to $\eta$  on $[0, T_*]$ with the rate $\sqrt{\s_n}$:
\bq\label{conv:root}
\|\eta_n-\eta\|_{L^\infty([0,T_*];H^{s-1})}+\|\eta_n-\eta\|_{L^2([0,T_*];H^{s-\mez})} \le  \sqrt{\s_n} \mathcal F(\|\eta_0\|_{H^s}, \ma^{-1}).
\eq
 If in addition $s\ge 2$, then we have the convergence with optimal rate $\s_n$:
\bq\label{conv:opt:intro}
\|\eta_n-\eta\|_{L^\infty([0,T_*];H^{s-2})}+\|\eta_n-\eta\|_{L^2([0,T_*];H^{s-\tdm})} \le {\s_n} \mathcal F(\|\eta_0\|_{H^s}, \ma^{-1}).
\eq
\end{theo}
The convergence \eqref{conv:root} holds for initial data in {\it any subcritical Sobolev spaces} $H^{1+\frac{d}{2}+}(\Rr^d)$. In particular, this allows for initial interfaces whose curvatures are {\it unbounded} in all dimensions and {\it not locally square integrable} in one dimension. The former is because $H(\eta_0)\in H^{-1+\frac{d}{2}+}(\Rr^d)\not\subset L^\infty(\Rr^d)$  and latter is due to the fact that in one dimension we have $H(\eta_0)\in H^{-\mez+\eps}(\Rr)\not\subset L^2_{loc}(\Rr)$. This appears to be the first result on vanishing surface tension  that can accommodate curvature singularities of the initial interface. On the other hand, the convergence \eqref{conv:opt:intro} has {\it optimal rate} $\s_n$ and holds under the additional condition that $s\ge 2$. This is only a condition in one dimension since $s>1+\frac{d}{2}\ge 2$ for $d\ge 2$. Note also that $s \ge 2$ is  the minimal regularity to ensure that the initial curvature is square integrable, yet it still allows for unbounded curvature.  See the technical remark \ref{rema:cd}. 

The proof of Theorem \ref{theo:main} exploits the Dirichlet-Neumann reformulation \cite{NgPa, Ng1} for the Muskat problem in a general setting. See also \cite{AlaMeuSme} for the one-fluid case. Part (i) of Theorem \ref{theo:main} is a uniform local well-posedness with repsect to surface tension. The key tool in proving this is paralinearization results for the Dirichlet-Neumann operator taken from \cite{NgPa, ABZ3}. The convergences \eqref{conv:root} and \eqref{conv:opt:intro} rely on contraction estimates for the Dirichlet-Neumann operator proved in \cite{NgPa} for a large Sobolev regularity range of the Dirichlet data. Together with \cite{NgPa} and \cite{Ng1}, Theorem \ref{theo:main}  provides a rather complete local regularity theory for (large) subcritical data.
\begin{rema}
In \cite{Ambrose2D} and \cite{Ambrose3D}, the first results on the zero surface tension limit for Muskat were obtained respectively in 2D and 3D for smooth initial data, i.e. $\Sigma_0\in H^{s_0}$  for some sufficiently large $s_0$. The interface is not necessarily a graph but if it is then the convergence estimates therein  translate into
\bq
\| \eta_n-\eta\|_{L^\infty([0, T_*]; H^1(\T^d))}\le C\sqrt{\s_n},\quad d=1, 2,
\eq 
which has the same rate as \eqref{conv:root}.
\end{rema}
\begin{rema}\label{rema:cd}
The condition $s\ge 2$ in \eqref{conv:opt:intro} is  due to the control of low frequencies in the paralinearization and contraction estimates for the Dirichlet-Neumann operator $G(\eta)f$ (see \eqref{def:DN-} for its definition). Precisely, the best currently available results (see subsections \ref{subsection:para} and \ref{subsection:contra} below) require $f\in  H^\sigma(\Rr^d)$ with $\sigma\ge \mez$. The proof of the  $L^\infty_t H^{s-2}_x$ convergence in \eqref{conv:opt:intro} appeals to these results with $\sigma=s-\tdm$. 
\end{rema}
\begin{rema}
It was proved in \cite{NgPa} that the Rayleigh-Taylor condition holds unconditionally in the following configurations:
\begin{itemize}
\item the one-phase problem without bottom or with Lipschitz bottoms;
\item the two-phase problem with constant viscosity ($\mu^+=\mu^-$).
\end{itemize}
When the Rayleigh-Taylor condition is violated, analytic solutions to the problem without surface tension exist \cite{SieCafHow}. The works \cite{SieTan, SieTanDai} and \cite{CenHou1, CenHou2} strongly indicate that these solutions are not limits  of solutions to the problem with surface tension. We also refer to \cite{GuoHalSpi} for the instability of the trivial solution (with surface tension) and to \cite{GanGarPatStr2} for the stability of  bubbles (without surface tension). 
\end{rema}
\begin{rema}
In Theorem  \ref{theo:main}, the initial data is fixed for all surface tension coefficients $\s_n$. In general, one can consider $\eta_n\vert_{t=0}=\eta_{n, 0}$ uniformly bounded bounded in $H^s(\Rr^d)$ such that the conditions in \eqref{cd:data:intro} hold uniformly in $n$. Then, for any $s>1+\frac{d}{2}$, we have 
\bq
\|\eta_n-\eta\|_{L^\infty([0,T_*];H^{s-1})}+\|\eta_n-\eta\|_{L^2([0,T_*];H^{s-\mez})} \le  \big(\sqrt{\s_n} +\| \eta_{n, 0}-\eta_0\|_{H^{s-1}}\big)\mathcal F(\|\eta_0\|_{H^s}, \ma^{-1}).
\eq
On the other hand, if $s>1+\frac{d}{2}$ and $s\ge 2$ then 
\bq
\|\eta_n-\eta\|_{L^\infty([0,T_*];H^{s-2})}+\|\eta_n-\eta\|_{L^2([0,T_*];H^{s-\tdm})} \le \big({\s_n}+\| \eta_{n, 0}-\eta_0\|_{H^{s-2}}\big) \mathcal F(\|\eta_0\|_{H^s}, \ma^{-1}).
\eq
\end{rema}
\begin{rema}
By interpolating the convergence estimate \eqref{conv:root} and the uniform bounds in \eqref{uniformbound:theo}, we obtain the vanishing surface tension limit in $H^{s'}$ for all $s'\in [s-1, s)$:
\bq
\|\eta_n-\eta\|_{L^\infty([0,T_*];H^{s'})}+\|\eta_n-\eta\|_{L^2([0,T_*];H^{s'+\mez})} \le  \s_n^{\frac{s-s'}{2}} \mathcal F(\|\eta_0\|_{H^s}, \ma^{-1}).
\eq
Convergence in the highest regularity $L^\infty_tH^s_x$ is more subtle and can possibly be established using  the  Bona-Smith type argument \cite{Bona}. This would imply in particular that $\eta$ is continuous in time with values in $H^s_x$, $\eta\in C_t H^s_x$. In the context of vanishing viscosity limit, this question was addressed in \cite{Masmoudi}, while convergence in lower Sobolev spaces (compared to initial data) was proved in \cite{Cons}. For gravity water waves, the Bona-Smith type argument was applied in \cite{Ng:flowmap} to  establish the continuity of the flow map in the highest regularity. 
\end{rema}
\begin{rema}
The proof of the local well-posedness in all subcritical Sobolev spaces  in \cite{NgPa} uses a parabolic regularization.  Theorem \ref{theo:main}  provides an alternative proof via regularization by vanishing surface tension. We stress that the assertions about $\eta$ in Theorem \ref{theo:main} do not make use of the local well-posedness results in \cite{NgPa}.
\end{rema}
The rest of this paper is organized as follows. In section \ref{section:prelim}, we recall the reformulation of the Muskat problem in terms of the Dirichlet-Neumann operator together with results on the Dirichlet-Neumann operator established in \cite{NgPa, ABZ3}. Section \ref{sectipn:unibound} is devoted to the proof of uniform-in-$\s$ {\it a priori} estimates. In section \ref{section:contractionJ}, we prove contraction estimates for the operators $J^\pm$ which arise in the reformulation of the two-phase problem. The proof of Theorem \ref{theo:main} is given in section \ref{section:proof}. Finally, in Appendix \ref{appendix}, we recall the symbolic paradifferential calculus and the G\r arding inequality for paradifferential operators. 

\section{Preliminaries}\label{section:prelim}
\subsection{Big-O notation}
If $X$ and $Y$ are Banach spaces and $T: X\to Y$ is a bounded linear operator, $T\in \mathcal{L}(X, Y)$,  with operator norm bounded by $A$, we write
\[
T=O_{X\to Y}(A).
\]
We define the space of operators of order $m \in \mathbb R$ in the scale of Sobolev spaces $H^s(\Rr^d)$:
\[
Op^m\equiv Op^m(\mathbb R^d) = \bigcap_{s \in \mathbb R} \mathcal L(H^s(\mathbb R^d),H^{s-m}(\mathbb R^d)).
\]
We  shall write $T=O_{Op^m}(A)$ when  $T\in Op^m$ and for all $s\in \Rr$, there exists $C=C(s)$ such that $T=O_{H^{s}\to H^{s-m}}(CA)$.
\subsection{Function spaces}
In the general setting, to reformulate the dynamics of the Muskat problem solely in terms of the interface, we require some function spaces. 

We shall always assume that $\eta\in W^{1, \infty}(\Rr^d)$ and either $\Gamma^\pm=\emptyset$ or $b^\pm \in \dot W^{1, \infty}(\Rr^d)$ with $\dist(\eta, b^\pm)>0$. Recall that the fluid domains $\Omega^\pm$ are given in \eqref{domain}. Define
\bq
\dot H^1(\Omega^\pm)=\{v\in L^1_{loc}(\Omega^\pm): \na_{x, y}v\in L^2(\Omega^\pm)\}/\Rr,\quad \| v\|_{\dot H^1(\Omega^\pm)}=\| \na_{x, y}v\|_{L^2(\Omega^\pm)}.
\eq
For any $\sigma \in \mathbb R$,  we define the `slightly homogeneous' Sobolev space
\bq
H^{1,\sigma}(\mathbb R^d) =\left\{ f\in L^1_{loc}(\Rr^d): \na f \in H^{\sigma-1}(\Rr^d)\right\}/\mathbb R,\quad \| f\|_{H^{1,\sigma}(\mathbb R^d)}=\| \na f\|_{H^{\sigma-1}(\mathbb R^d)}.
\eq
When $b^\pm\in \dot W^{1, \infty}(\Rr^d)$, we fix an arbitrary number $a\in (0, 1)$ and  define the `screened' fractional Sobolev spaces
\bq\label{def:wtHmez}
\wt  H^\mez_{\mp a(\eta-b^\pm)}(\Rr^d)=\left\{f\in L^1_{loc}(\Rr^d): \int_{\Rr^d}\int_{|k|\le \mp a(\eta-b^\pm)}\frac{|f(x+k)-f(x)|^2}{|k|^{d+1}}dkdx<\infty \right\}/\Rr.
\eq
According to Proposition 3.2 \cite{NgPa}, the spaces $\wt  H^\mez_{\mp a(\eta-b^\pm)}(\Rr^d)$ are independent of $\eta$ in $W^{1, \infty}(\Rr^d)$ satisfying $\dist(\eta, b^\pm)>h>0$. Thus, we can set
\bq
\wt H^\mez_\pm(\Rr^d)=\begin{cases}
\dot H^\mez(\Rr^d)\quad \text{if}~\Gamma^\pm=\emptyset,\\
\wt  H^\mez_{\mp a(\eta-b^\pm)}(\Rr^d)\quad\text{if}~b^\pm\in \dot W^{1, \infty}(\Rr^d).
\end{cases}
\eq
 It was proved in \cite{Stri, LeoTice} that  there exist  unique continuous trace operators
\bq\label{Tr}
\mathrm{Tr}_{\Omega^\pm\to \Sigma}: \dot H^1(\Omega^\pm)\to \wt H^\mez_\pm(\Sigma)\equiv \wt H^\mez_\pm(\Rr^d)
\eq
with norm depending only on $\| \eta\|_{\dot W^{1, \infty}(\Rr^d)}$ and $\| b^\pm\|_{\dot W^{1, \infty}(\Rr^d)}$.

Finally, for $\sigma> \mez$, we define
\bq
\wt H^\sigma_\pm(\Rr^d)=\wt H^\mez_\pm(\Rr^d)\cap H^{1, \sigma}(\Rr^d)
\eq
and equip it  with the norm $\| \cdot\|_{\wt H_\pm^\sigma}=\| \cdot\|_{\wt H^\mez_\pm}+\| \cdot\|_{H^{1, \sigma}}$.
\subsection{Dirichlet-Neumann formulation}
Given a function $f \in \mathbb R^d$, let $\phi$ solve the Laplace equation
\begin{equation}\label{phiSys}
\begin{cases}
\Delta_{x,y} \phi = 0 \text{ in } \Omega^-, \\ 
\phi = f \text { on } \Sigma, \\
\frac{\partial \phi}{\partial \nu^-} =0  \text{ on } \Gamma^-,
\end{cases}
\end{equation}  
with the final condition replaced by decay at infinity of $\phi$ if $\Gamma^- = \emptyset$.  Then, we define the (rescaled) Dirichlet-Neumann operator $G^- \equiv G^-(\eta)$ by 
\bq\label{def:DN-}
G^- f = \langle \nabla \eta\rangle \frac{\partial \phi^-}{\partial n}.
\eq
The operator $G^+(\eta)$ for the top fluid domain $\Omega^+$ is defined similarly. The solvability of \eqref{phiSys} is given in the next proposition. 
\begin{prop}[\protect{\cite{NgPa} Propositions 3.4 and 3.6}]\label{prop:variational}
Assume that either $\Gamma^-=\emptyset$ or $b^-\in  \dot W^{1, \infty}(\Rr^d)$. If $\eta\in  W^{1, \infty}(\Rr^d)$ and   $\dist(\eta, \Gamma^-)>h>0$, then for every $f\in \wt H^\mez_-(\Rr^d)$ there exists a unique variational solution $\phi^-\in \dot H^1(\Omega^-)$ to \eqref{phiSys}. Precisely, $\phi$ satisfies $\mathrm{Tr}_{\Omega^-\to \Sigma}\phi=f$,
\bq
\int_{\Omega^-}\na_{x, y}\phi\cdot\na_{x, y}\varphi dxdy=0\quad\forall \varphi\in \big\{v\in \dot H^1(\Omega^-):~\mathrm{Tr}_{\Omega^-\to \Sigma}v=0\big\},
\eq
together with the estimate 
\bq
\| \na_{x, y}\phi\|_{L^2(\Omega^-)}\le \cF(\| \na \eta\|_{L^\infty})\| f\|_{\wt H^\mez_-(\Rr^d)}
\eq
for some $\cF: \Rr^+\to \Rr^+$ depending only on $h$ and $\| \na_xb^-\|_{L^\infty(\Rr^d)}$.
\end{prop}
As the functions $b^\pm$ are fixed in $\dot W^{1, \infty}(\Rr^d)$, we shall omit the dependence on $\| \na_xb^-\|_{L^\infty(\Rr^d)}$ in various estimates below. 

The Muskat problem can be reformulated in terms of $G^\pm$ as follows.
\begin{prop}[\protect{\cite{Ng1} Proposition 1.1}]\label{prop:DNreform}

$(i)$ If $(u, p, \eta)$ solve the one-phase Muskat problem then $\eta:\Rr^d\to \Rr$  obeys the equation 
\bq\label{eq:eta1p}
\p_t\eta=-\frac{1}{\mu^-} G^{-}(\eta)\big(g\rho^-\eta+\s H(\eta)\big).
\eq
Conversely, if $\eta$  is a solution of \eqref{eq:eta1p} then the one-phase Muskat problem has a solution which admits $\eta$ as the free surface.

$(ii)$ If $(u^\pm, p^\pm, \eta)$ is a solution of the two-phase Muskat problem then 
\bq\label{eq:eta2p}
\p_t\eta=-\frac{1}{\mu^-}G^{-}(\eta)f^-,
\eq
where $f^\pm:=p^\pm\vert_{\Sigma}+\rho^\pm g \eta$ satisfy
\bq\label{system:fpm}
\begin{cases}
 f^--f^+=\g\eta+\s H(\eta),\\
\frac{1}{\mu^-}G^-(\eta)f^--\frac{1}{\mu^+}G^+(\eta)f^+=0.
\end{cases}
\eq
Conversely, if $\eta$ is a solution of  \eqref{eq:eta2p} where $f^\pm$ solve \eqref{system:fpm} then the two-phase Muskat problem has a solution which admits $\eta$ as the free interface.
\end{prop}
To make use of the results on the  Dirichlet-Neumann operator established in \cite{NgPa, Ng1, ABZ3}, it is convenient to introduce the linear operators
\bq\label{def:J}
 J ^\pm= J ^\pm(\eta):\quad v \mapsto f^\pm
 \eq
  where $f^\pm :\mathbb R^d \to \mathbb R$ are solutions to the system
\begin{equation}\label{JSystem}
\begin{pmatrix}
\text{Id} & -\mathrm{Id}\\
\mu^+ G^-(\eta) & -\mu^- G^+(\eta) 
\end{pmatrix}
\begin{pmatrix}
f^- \\ f^+
\end{pmatrix} = \begin{pmatrix}
v \\ 0
\end{pmatrix}.
\end{equation}
Introduce
 \bq\label{def:L}
 L = L(\eta) = \frac{\mu^+ + \mu^-}{\mu^-}G^- J ^-.
 \eq
For the two-phase case, $L$ coincides with $\frac{\mu^+ + \mu^-}{\mu^+}G^+  J ^+$ in view of \eqref{system:fpm}. Thus, writing $L = (\frac{\mu^-}{\mu^+ + \mu^-} +\frac{\mu^+}{\mu^+ + \mu^-})L$ yields the symmetric formula
\bq\label{form:L}
L = G^- J ^- + G^+ J ^+.
\eq
This formula holds for  the one phase problem  \eqref{eq:eta1p} as well. Indeed, when $\mu^+ = 0$ we have $ J ^+ = 0$ and $ J ^- = \text{Id}$, and hence $L = G^-$. In view of \eqref{def:L}, Proposition \ref{prop:DNreform} implies the following.
\begin{lemm} The  Muskat problem (both one-phase and two-phase) is equivalent to
\begin{equation}
 \partial_t\eta +\frac{1}{\mu^+ + \mu^-}L(\eta)( \mathfrak g \eta  + \mathfrak s H(\eta))= 0.\label{evol}
\end{equation}
\end{lemm} 
The next proposition gathers results on the existence and boundedness of the operators $J ^\pm$, $G^\pm$, and $L$ in Sobolev spaces.
\begin{prop}[\protect{\cite{ABZ3} Theorem 3.12, and \cite{NgPa} Propositions 3.8, 4.8, 4.10 and Remark 4.9}]\label{GBdd}
Let $\mu^+ \geq 0$ and $\mu^- > 0$. Assume $\dist(\eta, \Gamma^\pm)>h>0$.

(i) If $\eta \in W^{1, \infty}(\Rr^d)$ then 
\bq
\| J ^\pm\|_{\mathcal L(H^{\mez},\wt H^{\mez}_\pm)} + \|G^\pm\|_{\mathcal L(\wt H^{\mez}_\pm, H^{-\mez})}   + \|L\|_{\mathcal L(H^{\mez},H^{-\mez})}\le \mathcal F(\| \eta\|_{W^{1, \infty}}),
\eq
where $\mathcal F$ is nondecreasing  and depends only on $(h, \mu^\pm)$. 
 
(ii) If  $\eta \in H^s(\Rr^d)$ with $s > 1+\frac{d}{2}$ then for any $\sigma \in [\mez, s]$, we have 
\begin{align}\label{continuity:JGL}
\| J ^\pm\|_{\mathcal L(H^{\sigma},\wt H^{\sigma}_\pm)}+  \|G^\pm\|_{\mathcal L(\wt H^{\sigma}_\pm, H^{\sigma-1})}  +\|L\|_{\mathcal L(H^{\sigma},H^{\sigma-1})}  \leq \mathcal F(\|\eta\|_{H^s}),
\end{align}
where $\mathcal F$ is nondecreasing  and depends only on $(h, \mu^\pm, s, \sigma)$. 
\end{prop} 
\subsection{Paralinearization}\label{subsection:para}
Given a function $f$, define the operators
\begin{align}\label{Bop}
{\mathfrak B}^\pm f\equiv {\mathfrak B}^\pm(\eta)f & = \langle \nabla \eta\rangle^{-2}(\nabla \eta \cdot \nabla  + G^\pm(\eta))f, \\ \label{Vop}
{\mathfrak V}^\pm f \equiv {\mathfrak V}^\pm(\eta)f & =(\nabla - \nabla \eta {\mathfrak B}^\pm )f.
\end{align} 
We note here that ${\mathfrak B}^\pm f = \partial_y \phi^\pm |_\Sigma$ and ${\mathfrak V}^\pm f= \nabla_x \phi^\pm |_{\Sigma}$, where $\phi$ solves \eqref{phiSys}. Moreover, as a consequence of \eqref{continuity:JGL} and product rules, we have 
\bq\label{continuity:BV}
\| {\mathfrak B}^\pm(\eta)\|_{\wt H^\sigma_\pm\to H^{\sigma-1}}+\| {\mathfrak V}^\pm(\eta)\|_{\wt H^\sigma_\pm\to H^{\sigma-1}}\le \cF(\| \eta\|_{H^s}) 
\eq
for all $s>1+\frac d2$ and $\sigma\in [\mez, s]$.

 The principal symbol of $G^-(\eta)$ is  
\begin{equation}\label{ld}
\lambda(x, \xi) = \sqrt{\langle \nabla \eta\rangle^2 |\xi|^2  - (\nabla \eta \cdot \xi)^2},
\end{equation}
while that of $G^+(\eta)$ is $-\ld(x, \xi)$. Note that $\lambda(x, \xi) \ge |\xi|$ with equality when $d=1$. Next we record results on paralinearization of $G^\pm(\eta)$.
\begin{theo}[\protect{\cite{ABZ3} Propostion 3.13,  \cite{NgPa} Theorem 3.18}] \label{DNpara}
Let $s>1+\frac{d}{2}$ and  let $\delta \in (0,\mez]$ satisfy $\delta<s-1-\frac d2$.   Assume that $\eta\in H^s$ and $\dist(\eta, \Gamma^\pm)>h>0$.

(i)  For any $\sigma \in [\mez, s-\delta]$, there exists a nondecreasing function $\mathcal F$ depending only on $(h, s,\sigma,\delta)$ such that
\begin{align}
&\mp G^\pm(\eta) = T_\lambda + O_{\wt H^\sigma_\pm \to  H^{\sigma - 1 + \delta}}(\mathcal F(\|\eta\|_{H^s}))\label{DNprincipal},\\
& L(\eta)  = T_\lambda +O_{H^\sigma\to H^{\sigma - 1+\delta}}( \mathcal F(\|\eta\|_{H^s})). \label{Lprincipal}
\end{align}

(ii) For any $\sigma \in [\mez, s]$, there exists a nondecreasing function $\mathcal F$ depending only on $(h, s,\sigma,\delta)$ such that for all $f \in \tilde  H^\sigma_\pm(\Rr^d)$, we have
\begin{align}
&\mp G^\pm(\eta)f = T_\lambda (f - T_{{\mathfrak B}^\pm f}\eta) - T_{{\mathfrak V}^\pm f}\cdot\nabla \eta + O_{\wt H^\sigma_\pm \to H^{\sigma-\mez}}(\mathcal F(\|\eta\|_{H^s})(1 + \| \eta\|_{H^{s+ \mez - \delta}}))f\label{DNgoodvar},\\
& L(\eta)f = T_\lambda(f - T_{\lbb {\mathfrak B} J  \rbb f}\eta) - T_{\lbb {\mathfrak V} J  \rbb f}\cdot \nabla \eta  + O_{H^\sigma \to H^{\sigma - \mez}}(\mathcal F(\|\eta\|_{H^s})(1 + \| \eta\|_{H^{s + \mez-\delta}}))f.\label{Lgoodvar}
\end{align}
\end{theo}
\begin{proof}
 \eqref{DNprincipal}  was proved in Proposition  3.13 in \cite{ABZ3} for $\sigma\in [\mez, s-\mez]$ but its proof allows for $\sigma\in [\mez, s-\delta]$. On the other hand,  \eqref{DNgoodvar} was proved in Theorem 3.17  in \cite{NgPa}. Let us prove \eqref{Lprincipal} and \eqref{Lgoodvar}. We recall from \eqref{form:L} that $L= G^- J ^-+G^+J^+$. Applying \eqref{DNprincipal} we have
\[
\mp G^\pm J^\pm f = T_\lambda J^\pm f+ O_{\wt H^\sigma_\pm \to  H^{\sigma - 1 + \delta}}(\mathcal F(\|\eta\|_{H^s}))J^\pm f.
\]
By virtue of \eqref{continuity:JGL} we have $\| J ^\pm\|_{\mathcal L(H^{\sigma},\wt H^{\sigma}_\pm)}\le \cF(\| \eta\|_{H^s})$, and thus
\[
\mp G^\pm J^\pm f = T_\lambda J^\pm f+ O_{ H^\sigma \to  H^{\sigma - 1 + \delta}}(\mathcal F(\|\eta\|_{H^s})) f.
\]
It follows that
\begin{align*}
L(\eta)f&=T_\lambda \lbb J\rbb f + O_{ H^\sigma \to  H^{\sigma - 1 + \delta}}(\mathcal F(\|\eta\|_{H^s}))f\\
&=T_\lambda f +  O_{ H^\sigma \to  H^{\sigma - 1 + \delta}}(\mathcal F(\|\eta\|_{H^s}))f,
\end{align*}
where in the second equality we have used  the fact that $ \lbb J\rbb=\mathrm{Id}$. This completes the proof of \eqref{Lprincipal}. Finally, \eqref{Lgoodvar} can be proved similarly upon using the paralinearizaion \eqref{DNgoodvar}.
\end{proof}
Finally, the mean curvature operator $H(\cdot)$, defined by \eqref{def:H}, can be paralinearized as follows.   
\begin{prop}[\protect{\cite{Ng1} Proposition 3.1}]\label{Hpara}
Let $s > 1 + \frac{d}{2}$ and let $\delta \in (0,\mez]$ satisfying $\delta<s-1-\frac d2$.  Then there exists a nondecreasing function $\cF$ depending only on $s$ such that 
\begin{equation}\label{lin:H}
H(\eta)  = T_{l}\eta +  O_{H^{s+\tdm} \to H^{s-\mez + \delta}}(\mathcal F(\|\eta\|_{H^s}) )\eta
\end{equation}
where
\bq
l = \langle \nabla \eta\rangle^{-3}\lambda^2.
\eq
In addition, if $\sigma\ge -1$ then 
\begin{equation}\label{boundH}
\|H(\eta)\|_{H^\sigma} \leq \mathcal F(\|\eta\|_{H^s})\|\eta\|_{H^{\sigma+2}}.
\end{equation}
\end{prop}
\subsection{Contraction estimates}\label{subsection:contra}
Let $s>1+\frac d2$ and consider  $\eta_j\in H^s(\Rr^d)$ satisfying $\dist(\eta_j, \Gamma^\pm)>h>0$, $j=1, 2$. We have the following contraction estimates for $G^\pm(\eta_1)-G^\pm(\eta_2)$.
\begin{theo}[\protect{\cite{NgPa} Corollary 3.25 and Proposition 3.31}]\label{theo:contractionDN}
For any $\sigma\in [\mez, s]$, there exists a nondecreasing function $\cF:\Rr^+\to \Rr^+$ depending only on $(h, s, \sigma)$ such that
\bq
\label{contraction:DN2}
\| G^\pm(\eta_1)-G^\pm(\eta_2) \|_{\wt H^s_\pm\to H^{\sigma-1}}\le \cF\big(\|(\eta_1, \eta_2)\|_{H^s}\big)\| \eta_1-\eta_2\|_{H^\sigma}
\eq
and
\bq\label{contraction:DN3}
\| G^\pm(\eta_1)-G^\pm(\eta_2)  \|_{\wt H^\sigma_\pm\to H^{\sigma-1}}\le 
\mathcal{F}\big(\| (\eta_1, \eta_2) \|_{H^s}\big)\| \eta_1-\eta_2\|_{H^s}.
\eq
\end{theo}
\begin{theo}[\protect{\cite{NgPa} Theorem 3.24}]\label{DNcont}
Let $\delta\in (0, \mez]$ satisfy $\delta<s-1-\frac d2$. Let $\sigma \in [\mez+\delta, s]$. For any $f \in \wt H^s_\pm$, there exists a nondecreasing function $\mathcal F$ depending only on $(h, s, \sigma)$ such that
\begin{equation}
\begin{aligned}
\mp  \big(G^\pm(\eta_1)f-G^\pm(\eta_2) f\big) &= -T_{\lambda_1 \mathfrak B^\pm(\eta_1) f}(\eta_1-\eta_2) - T_{\mathfrak V^\pm(\eta_1) f}\cdot\nabla (\eta_1-\eta_2)\\
&\qquad+O_{\wt H^s_\pm \to H^{\sigma-1}}\Big(\mathcal \cF\big(\|(\eta_1, \eta_2)\|_{H^s}\big)\|\eta_1-\eta_2\|_{H^{\sigma-\delta}}\Big)f,
\end{aligned}
\end{equation}
where $\ld_1$ is defined by \eqref{ld} with $\eta=\eta_1$.
\end{theo}
\section{Uniform {\it a priori} estimates}\label{sectipn:unibound}
Conclusion (i) in Theorem \ref{theo:main} concerns the uniform local well-posedness of the Muskat problem with surface tension. Key to that is the following  {\it a priori} estimates that are uniform in the vanishing surface tension limit $\s\to 0$. 
\begin{prop}\label{apriori}
Let $s > 1+\frac{d}{2}$, $\mu^- > 0$, $\mu^+ \geq 0$, $\mathfrak s > 0$, and $h>0$. Suppose 
\bq\label{Regcond}
\eta\in C([0, T]; H^s(\Rr^d))\cap L^2([0, T]; H^{s+\tdm}(\Rr^d))
\eq
 is a solution to $\eqref{evol}$ with initial data $\eta_0\in H^s(\Rr^d)$ such that 
\begin{align}\label{RTcond}
& \inf_{t\in [0, T]}\inf_{x\in \Rr^d} \mathrm{RT}(\eta(t)) > \ma>0,\\ \label{Sepcond}
&\inf_{t\in [0, T]}\dist(\eta(t), \Gamma^\pm)> h.
\end{align}
Then, there exists a nondecreasing function $\cF:\Rr^+\times \Rr^+\to \Rr^+$  depending only on $(h, s, \mu^\pm)$ such that
\begin{equation}\label{uni:Hs}
\begin{aligned}
&\|\eta\|_{L^\infty([0,T];H^s)}  \leq \|\eta_0\|_{H^s} \exp\Big((\mathfrak s  + \mathfrak g)T\mathcal F\big(\|\eta\|_{L^\infty([0,T];H^s)},\mathfrak a^{-1}\big)\Big)
\end{aligned}
\end{equation}
and 
\bq\label{uni:dissipation}
\mathfrak s \|\eta\|_{L^2([0,T];H^{s+\tdm})}^2  +\mathfrak g \|\eta\|_{L^2([0,T];H^{s+\mez})}^2\le \mathcal F_1\left( \|\eta_0\|_{H^s} \exp\Big((\mathfrak s  + \mathfrak g)T\mathcal F\big(\|\eta\|_{L^\infty([0,T];H^s)},\mathfrak a^{-1}\big)\Big),\mathfrak a^{-1}\right)
\eq
where $\cF_1(m, n)=m^2\cF(m, n)$.
\end{prop}
\begin{proof}
Set $ B= \lbb \mathfrak B(\eta)  J (\eta)\rbb \eta$ and $V= \lbb \mathfrak V(\eta)  J (\eta)\rbb \eta$. We shall write $\mathcal{Q}=\mathcal{Q}(t)=\cF(\| \eta(t)\|_{H^s}, \ma^{-1})$ when $\cF(\cdot, \cdot)$ is nondecreasing and depends only on $(h, s, \mu^\pm)$. Note that $\cF$ may change from line to line. From \eqref{evol} we have
\begin{equation}\label{eq:eta:2}
(\mu^+ + \mu^-)\partial_t\eta +\g L(\eta)\eta  + \s L(\eta)H(\eta)= 0
\eq
for both the one-phase and two-phase problems. Fix  $\delta\in \big(0, \min(\mez, s-1-\frac d2)\big)$.  
By virtue of the paralinearization \eqref{Lgoodvar} (with $\sigma=s$) we have
\[
\g L(\eta)\eta=\mathfrak g  (T_\lambda(\eta - T_{B}\eta) -T_V\cdot \nabla \eta)+O_{H^s \to H^{s-\mez}}\big(\mathfrak g \mathcal Q(1+\|\eta\|_{H^{s+\mez-\delta}})\big)\eta.
\]
On the other hand, \eqref{Lprincipal} (with $\sigma=s-\mez$) together with \eqref{boundH} gives
\[
 \s L(\eta)H(\eta)=\s T_\lambda H(\eta)+O_{H^{s+\tdm}\to H^{s-\tdm+\delta}}(\s\mathcal F(\|\eta\|_{H^s}))\eta
 \]
 Combining this with the linearization \eqref{lin:H} for $H(\eta)$ yields
 \[
 \s L(\eta)H(\eta)=\s T_\lambda T_l\eta+O_{H^{s+\tdm}\to H^{s-\tdm+\delta}}(\s\mathcal F(\|\eta\|_{H^s}))\eta,
 \]
 where we have applied Theorem \ref{theo:Op} to have $T_\ld=O_{Op^1}(\cF(\|\eta\|_{H^s}))$.
 
 Then in view of \eqref{eq:eta:2} we obtain
\begin{equation}
\begin{aligned}
&(\mu^+ + \mu^-)\partial_t \eta + \mathfrak s T_\lambda T_{l} \eta +  \mathfrak g  (T_\lambda(\eta - T_{B}\eta) -T_V\cdot \nabla \eta) \\
& \ \ \  \ \ \ \ \ \ \  = O_{H^s \to H^{s-\mez}}\Big(\mathfrak g \mathcal F(\|\eta\|_{H^s})  (1+\|\eta_s\|_{H^{\mez-\delta}})\Big)\eta +  O_{H^{s+\tdm} \to H^{s-\tdm+\delta}}(\mathfrak s \mathcal F(\|\eta\|_{H^s}))\eta. \label{para1}
\end{aligned}
\end{equation}
Note that $\ld\in \Gamma^1_\delta$, $l\in \Gamma^2_\delta$ and $(B, V)\in W^{1+\delta, \infty}\subset\Gamma^0_\delta$, with seminorms bounded by $\cF(\| \eta\|_{H^s})$.  The symbolic calculus in Theorem \ref{theo:sc} then gives
\begin{align*}
&T_\lambda T_l \eta = T_{\lambda l}\eta  
 + O_{Op^{3-\delta}}(\mathcal F(\|\eta\|_{H^s})),\\
 &T_{\lambda}(\text{Id} - T_{ B}) =  T_{\lambda(1-B)} + O_{Op^{1-\delta}}(\mathcal F(\|\eta\|_{H^s})),\\
 &T_{V}\cdot \nabla=iT_{\xi\cdot V}=  i\mathrm {Re}(T_{\xi \cdot V}) +O_{Op^{1-\delta}}(\mathcal F(\|\eta\|_{H^s})),
\end{align*}
where  $\mathrm{Re}(T_{\xi\cdot V}) = \frac{1}{2}(T_{\xi\cdot V} + T_{\xi\cdot V}^*)$. It then follows from \eqref{para1} that
\begin{equation}
\begin{aligned}
&(\mu^+ + \mu^-)\partial_t \eta + \mathfrak s T_{\lambda l} \eta +  \mathfrak g \Big(T_{\lambda(1-B)}\eta -i\text{Re}(T_{\xi\cdot V}) \eta\Big) \\
& \ \ \  \ \ \ \ \ \ \  = O_{H^s \to H^{s-\mez}}(\mathfrak g \mathcal F(\|\eta\|_{H^s})  (1+\|\eta_s\|_{H^{\mez-\delta}}))\eta +  O_{H^{s+\tdm} \to H^{s-\tdm+\delta}}(\mathfrak s \mathcal F(\|\eta\|_{H^s}))\eta\label{para2}.
\end{aligned}
\end{equation}
We set $\eta_s = \langle D\rangle^s \eta$.  Appealing to Theorem \ref{theo:sc} again we have
\begin{align*}
&[\langle D\rangle^s,T_{\lambda l}] = O_{Op^{s+3-\delta}}(\mathcal F(\|\eta\|_{H^s})),\\
&[\langle D\rangle^s,T_{\lambda(1-B)}] =   O_{Op^{s+1-\delta}}(\mathcal F(\|\eta\|_{H^s}, \ma^{-1})), \\
&[\langle D\rangle^s,\mathrm{Re}(T_{\xi \cdot V})]=    O_{Op^{s+1-\delta}}(\mathcal F(\|\eta\|_{H^s})),
\end{align*}
where $[A, B]=AB-BA$ and in the second line we have used the lower bound \eqref{RTcond} for $(1-B)$ together with the fact that $\ld\ge |\xi|$. Note that we have adopted the convention that $\mathcal F(\|\eta\|_{H^s})\equiv \mathcal F(\|\eta\|_{H^s}, 0)$. This implies
\begin{equation}
\begin{aligned}
&(\mu^+ + \mu^-)\partial_t \eta_s + \mathfrak s T_{\lambda l} \eta_s +  \mathfrak g \Big(T_{\lambda(1-B)}\eta_s -\mathrm{Re}(T_{\xi\cdot V})\eta_s\Big) \\
& \ \ \  \ \ \ \ \ \ \  =  O_{L^2 \to H^{-\mez}}(\mathfrak g \mathcal Q(1+\|\eta_s\|_{H^{\mez-\delta}})) \eta_s+  O_{Op^{1-\delta}}(\mathfrak g \mathcal Q)\eta_s+ O_{H^{\tdm} \to H^{-\tdm+\delta}}(\mathfrak s \mathcal Q)\eta_s \label{para3}.
\end{aligned}
\end{equation}
Since $i\text{Re}(T_{\xi\cdot V})$ is skew-adjoint, by testing \eqref{para3} against $\eta_s$, we obtain
\begin{equation}
\begin{aligned}
&\frac{(\mu^+ + \mu^-)}{2}\frac{d}{dt}\| \eta_s\|_{L^2} ^2+ \mathfrak s(T_{{\lambda l}}\eta_s,\eta_s)_{L^2}+  \mathfrak g(T_{{\lambda(1-B)}}\eta_s,\eta_s)_{L^2} \\
&\ \ \  \ \ \ \ \ \ \ \leq \mathcal Q\Big\{\mathfrak g\Big[(1+\|\eta_s\|_{H^{\mez-\delta}}) \|\eta_s\|_{L^2}\|\eta_s\|_{H^{\mez}} + \|\eta_s\|_{H^{\mez-\delta}}\|\eta_s\|_{H^{\mez}}\Big]+ \mathfrak s \|\eta_s\|_{H^{\tdm-\delta}}\|\eta_s\|_{H^{\tdm}}\Big\},
\end{aligned}\label{energy1}
\end{equation}
where $(\cdot, \cdot)_{L^2}$ denotes the $L^2$ pairing.  The term involving $1+ \|\eta_s\|_{H^{\mez-\delta}}$ is treated as follows:
\begin{equation}
\begin{aligned}
(1+ \|\eta_s\|_{H^{\mez-\delta}}) \|\eta_s\|_{L^2}\|\eta_s\|_{H^{\mez}} &\leq \|\eta_s\|_{L^2}\|\eta_s\|_{H^{\mez}}+ \|\eta_s\|_{L^2}\|\eta_s\|_{H^{\mez-\delta}}\|\eta_s\|_{H^{\mez}} \\
&\leq \mathcal Q\|\eta_s\|_{H^{\mez-\delta}}\|\eta_s\|_{H^{\mez}}.
\end{aligned}
\end{equation}
In view of \eqref{RTcond} and the fact that $\ld(x, \xi)\ge |\xi|$, we have the lower bounds
\[
\lambda(1-B)\ge \ma|\xi|,\quad l\lambda=\langle \na \eta\rangle^{-3}\ld^3\ge \frac{1}{\langle \|\eta\|_{W^{1, \infty}}\rangle^3}|\xi|^3.
\]
Moreover, $\lambda(1-B)\in \Gamma^1_\delta$ and $l\lambda\in \Gamma^3_\delta$ with seminorms bounded by $\cF(\| \eta\|_{H^s})$. Then applying  the G\r arding's inequality \eqref{wgi} gives
\begin{align}
\|\Psi(D)\eta_s\|_{H^{\tdm}}^2&\leq \mathcal Q\Big(\|\eta_s\|_{H^\tdm}\|\eta_s\|_{H^{\tdm-\delta}}+  (T_{l\lambda}\eta_s,\eta_s)_{L^2}\Big),\\
\|\Psi(D)\eta_s\|_{H^{\mez}}^2 &\leq \mathcal Q \Big(\|\eta_s\|_{H^\mez}\|\eta_s\|_{H^{\mez-\delta}}^2 + (T_{\lambda(1-B)}\eta_s,\eta_s)_{L^2}\Big),
\end{align}
where $\Psi(D)$ denotes the Fourier multiplier with symbol $\Psi$ defined by \eqref{cond.psi}. In addition, we have
\[
\| u\|_{H^r}\le C(\| \Psi(D)u\|_{H^r}+\| u\|_{L^2})\quad\forall r\in \Rr.
\]
Thus, \eqref{energy1} amounts to
\begin{equation}
\begin{aligned}
&\frac{(\mu^+ + \mu^-)}{2}\frac{d}{dt}\| \eta_s\|_{L^2} ^2+ \frac{1}{\mathcal Q}\Big({\mathfrak s}\| \eta_s\|_{H^{\tdm}}^2 +  {\mathfrak g}\| \eta_s\|_{H^{\mez}}^2\Big) \\
& \ \ \ \ \ \ \ \ \ \ \leq \mathcal Q(\mathfrak g\|\eta_s\|_{H^{\mez-\delta}}\|\eta_s\|_{H^{\mez}}+  \mathfrak s \|\eta_s\|_{H^{\tdm-\delta}}\|\eta_s\|_{H^{\tdm}}).
\end{aligned}\label{energy2}
\end{equation}
We use Young's inequality and interpolation as follows:
\begin{align}
\|\eta_s\|_{H^{\mez-\delta}}\|\eta_s\|_{H^{\mez}} &\leq \|\eta_s\|_{L^2}^{2\delta}\|\eta_s\|_{H^{\mez}}^{2(1-\delta)} \leq (10\mathcal Q)^{\frac{2(1-\delta)}{\delta}}\|\eta_s\|_{L^2}^2 + \frac{1}{100\mathcal Q^2}\|\eta_s\|_{H^{\mez}}^2
\end{align}
and similarly,
\begin{align}
\|\eta_s\|_{H^{\tdm-\delta}}\|\eta_s\|_{H^{\tdm}}\le \|\eta_s\|_{L^2}^{2\frac{\delta}{3}}\|\eta_s\|_{H^{\tdm}}^{2(1-\frac{\delta}{3})}\leq (10\mathcal Q)^{\frac{2(3-\delta)}{\delta}}\|\eta_s\|_{L^2}^2 + \frac{1}{100\mathcal Q^2}\|\eta_s\|_{H^{\tdm}}^2.
\end{align}
Applying these inequalities to \eqref{energy2}, and then subtracting terms involving $\|\eta_s\|_{H^{\mez}}$, we obtain for  a larger $\mathcal{Q}$ if needed that
\begin{equation}
\frac{\mu^+ + \mu^-}{2}\frac{d}{dt} \|\eta_s\|_{L^2}^2 + \frac{1}{\mathcal Q}(\mathfrak s \| \eta_s\|_{H^{\tdm}}^2 + \mathfrak g \| \eta_s\|_{H^{\mez}}^2) \leq (\mathfrak g + \mathfrak s) \mathcal Q\|\eta_s\|_{L^2}^2.\label{energy3}
\end{equation}
A Gr\"onwall's argument then leads to 
\[
\begin{aligned}
&\|\eta\|_{L^\infty([0,T];H^s)}^2 + \frac{1}{\mathcal{Q}_T}\Big(\mathfrak s \|\eta\|_{L^2([0,T];H^{s+\tdm})}^2  +\mathfrak g \|\eta\|_{L^2([0,T];H^{s+\mez})}^2\Big) \leq \|\eta_0\|_{H^s}^2 \exp\Big((\mathfrak s  + \mathfrak g)T\mathcal Q_T\Big),
\end{aligned}
\]
where $\mathcal Q_T=\mathcal F(\|\eta\|_{L^\infty([0,T];H^s)},\mathfrak a^{-1})$ with $\cF$ depending only on $(h, s, \mu^\pm)$. In particular, we have the $H^s$ estimate \eqref{uni:Hs}. As for the dissipation estimate, we have 
\begin{align*}
&\mathfrak s \|\eta\|_{L^2([0,T];H^{s+\tdm})}^2  +\mathfrak g \|\eta\|_{L^2([0,T];H^{s+\mez})}^2\le \|\eta_0\|_{H^s}^2 \exp\Big((\mathfrak s  + \mathfrak g)T\mathcal F\big(\|\eta\|_{L^\infty([0,T];H^s)},\mathfrak a^{-1}\big)\Big)\mathcal{Q}_T.
\end{align*}
On the other hand, plugging \eqref{uni:Hs} into $\mathcal{Q}_T$ gives
\begin{align*}
\mathcal{Q}_T&=\mathcal F(\|\eta\|_{L^\infty([0,T];H^s)},\mathfrak a^{-1})\le \mathcal F\left( \|\eta_0\|_{H^s} \exp\Big((\mathfrak s  + \mathfrak g)T\mathcal F\big(\|\eta\|_{L^\infty([0,T];H^s)},\mathfrak a^{-1}\big)\Big),\mathfrak a^{-1}\right)
\end{align*}
Therefore, upon setting $\cF_1(m, n)=m^2\cF(m, n)$ we obtain 
\begin{align*}
\mathfrak s \|\eta\|_{L^2([0,T];H^{s+\tdm})}^2  +\mathfrak g \|\eta\|_{L^2([0,T];H^{s+\mez})}^2\le \mathcal F_1\left( \|\eta_0\|_{H^s} \exp\Big((\mathfrak s  + \mathfrak g)T\mathcal F\big(\|\eta\|_{L^\infty([0,T];H^s)},\mathfrak a^{-1}\big)\Big),\mathfrak a^{-1}\right)
\end{align*}
which finishes the proof of \eqref{uni:dissipation}.
\end{proof}

\section{Contraction  estimates for $J^\pm$}\label{section:contractionJ}
Our goal in this section is to prove  contraction estimates for $J^\pm(\eta)$ at two different  surfaces $\eta_1$ and $\eta_2$. This is only a question for the two-phase problem since for the one-phase problem we have $J^-=\mathrm{Id}$ and $J^+\equiv 0$.  Given an object $X$ depending  on $\eta$, we shall denote
 $X_j = X|_{\eta = \eta_j}$ and the difference
 \[
 X_\delta = X_1 - X_2.
  \]
\begin{prop}\label{gamContLem} 
Let $s > 1 + \frac{d}{2}$ and consider  $\eta_j\in H^s(\Rr^d)$ satisfying $\dist(\eta_j, \Gamma^\pm)> h>0$, $j=1, 2$.  For any $\sigma\in [\mez, s]$,  there exists $\cF:\Rr^+\to \Rr^+$ depending only on $(h, s, \sigma, \mu^\pm)$ such that
\begin{align}
&\| J ^\pm_\delta \|_{ H^s\to \wt H^\sigma_\pm} \leq \mathcal F(N_s)\|\eta_\delta\|_{H^\sigma},\label{gamCont12}\\
&\| J ^\pm_\delta \|_{ H^\sigma\to \wt H^\sigma_\pm} \leq \mathcal F(N_s)\|\eta_\delta\|_{H^s}\label{gamCont11},
\end{align} 
where we denoted 
\bq\label{Ns}
 N_s = \|(\eta_1, \eta_2)\|_{H^s}.
 \eq 
\end{prop}
We shall prove Proposition \ref{gamContLem} for the most general case of two fluids and with bottoms, i.e.  $\mu^+>0$ and $\Gamma^\pm\ne \emptyset$. Adaption to the other cases is straightforward.   
\subsection{Flattening the domain}
 There exist $\eta_*^\pm\in C_b^{s+100}(\Rr^d)$ such that
\bq\label{def:eta*}
b^-(x)+\frac{h}{2}\le \eta^-_*(x)\le \eta_j(x)-\frac{h}{2},\quad \eta_j(x)\le \eta_*^+(x)-\frac{h}{2}\le b^+(x)-h\quad\forall x\in \Rr^d
\eq
and for some $C=C(h, s, d)$,
\bq\label{norm:eta*}
\| \eta_*^\pm\|_{C_b^{s+100}(\Rr^d)}\le C(1+\| \eta_1\|_{L^\infty}+\| \eta_2\|_{L^\infty}).
\eq
For $j=1, 2$ we set $\Omega_*=\Omega^+_{j, *}\cup \Omega^-_{j, *}$ where
 \begin{equation}\label{divideOmegaj}
\Omega_{*, j}^\pm = \{(x,y):x\in \Rr^d,  \pm\eta(x)\le\pm y\le\pm\eta^\pm_*(x)\}.
 \end{equation}
Note that $\Omega_*=\{(x,y):x\in \Rr^d,  \eta^-_*(x)\le y\le \eta^+_*(x)\}$ is independent of $j\in \{1, 2\}$. For small $\tau>0$ to be chosen, define $\rho_j(x, z): \Rr^d\times [-1, 1]$ by
 \begin{equation}\label{diffeo}
 \varrho_j(x,z)=
  (1-z^2)e^{-\tau |z|\langle D_x \rangle }\eta_j(x) -\mez z(1-z) \eta^-_*(x)+\mez z(1+z)\eta^+_*(x).
\end{equation}
%
%
\begin{lemm}\label{diffeos}
There exists $K>0$ depending only on $(s, d)$ such that if  $\tau K\| \eta_j\|_{H^s}\le \frac{h}{12}$ then 
  \bq\label{lower:varrho}
  \p_z\varrho_j(x, z)\ge \frac{h}{12}\quad a.e.~(x, z)\in \Rr^d\times (-1, 1),~j=1, 2.
  \eq
For $j =1,2$, the mapping
\[
\Phi_j: \Rr^d\times [-1, 1] \ni (x, z) \mapsto (x, \varrho_j(x, z))\in \Omega_*
\]
is a Lipschitz diffeomorphism and  respectively maps $\Rr^d\times [0, 1]$ and $\Rr^d\times [-1, 0]$ onto $\Omega_{*, j}^+$ and  $\Omega_{*, j}^-$. Moreover, there exists $C=C(h, s, d)$ such that
\begin{align}\label{est:Phij}
&\| \na \Phi_j\|_{L^\infty(\Rr^d\times (-1, 1))}\le C(1+N_s),\\  \label{est:Phi-1}
&\| \na(\Phi_j^{-1})\|_{L^\infty(\Omega_*)}\le C(1+N_s).
\end{align}
\end{lemm}
\begin{proof}
We first note that $\Phi_j(x, 0)=\Sigma_j=\{(x, \eta_j(x)): x\in \Rr^d\}$, $\Phi_j(x, 1)=\{(x, \eta^+_*(x)): x\in \Rr^d\}$ and $\Phi_j(x, -1)=\{(x, \eta^-_*(x)): x\in \Rr^d\}$. Thus, in order to prove that $\Phi_j$ is one-to-one and onto, it suffices to prove that $\p_z\varrho_j(x, z)\ge c>0$ for a.e. $(x, z)\in \Rr^d\times (-1, 1)$. For $z\in (-1, 1)\setminus\{0\}$ we have
\begin{align*}
\p_z\varrho_j(x, z)&=\mez(1-2z)(\eta_j(x)-\eta^-_*(x))+\mez(1+2z)(\eta^+_*(x)-\eta_j(x))\\
&\quad-2z(e^{-\tau |z|\langle D_x \rangle }-1)\eta_j(x)-\text{sign}(z)\tau(1-z^2)e^{-\tau |z|\langle D_x \rangle }\langle D_x \rangle\eta_j(x).
\end{align*}
For $z\in [\frac13, 1]$, $1-2z\in [-1, \frac13]$ and $1+2z\in [\frac53, 3]$. In addition, by \eqref{def:eta*} we have $\mp(\eta_-\eta_*^\pm)\ge h/2$. Consequently, 
\bq\label{dzvarrho:0}
\mez(1-2z)(\eta_j(x)-\eta^-_*(x))+\mez(1+2z)(\eta^+_*(x)-\eta_j(x))\ge\frac{h}{6}.
\eq
Similarly we obtain \eqref{dzvarrho:0} for $z\in (-1, -\frac13)$ and $z\in [-\frac13, \frac13]\setminus\{0\}$. Next writing 
\[
(e^{-\tau |z|\langle D_x \rangle }-1)\eta_j(x)=-\tau \int_0^{|z|}e^{-\tau z'\langle D_x \rangle }\langle D_x \rangle \eta_j(x)dz'
\]
we obtain that 
\[
\p_z\varrho_j(x, z)\ge \frac{h}{6}-\tau K\| \eta_j\|_{H^s}\quad\forall (x, z)\in\Rr^d\times (-1, 1)
\]
for some constant $K=K(s, d)$. Note that   the condition $s>1+\frac d2$ has been used. Choosing $\tau>0$ such that $\tau K\| \eta_j\|_{H^s}\le h/12$  gives $\p_z\varrho_j(x, z)\ge h/12$ for a.e. $(x, z)\in \Rr^d\times (-1, 1)$.

Since 
\begin{align*}
\p_z\varrho_j(x, z)&=-2ze^{-\tau |z|\langle D_x \rangle }\eta_j(x)-\text{sign}(z)\tau(1-z^2)e^{-\tau |z|\langle D_x \rangle }\langle D_x \rangle\eta_j(x)\\
&-\mez(1-2z)\eta^-_*(x)+\mez(1+2z)\eta^+_*(x),
\end{align*}
there exists $K'=K'(s, d)$ such that
\[
\| \p_z\varrho_j(x, z)\|_{W^{1, \infty}(\Rr^d\times (-1, 1))}\le K'\| \eta_j\|_{H^{s-1}}+K'\tau \|\eta_j\|_{H^s}+K'(1+N_{s-1}),
\]
where we have used \eqref{norm:eta*}. Then in view of the fact that $\tau \|\eta_j\|_{H^s}\le hK^{-1}$, we obtain 
\bq\label{est:varrho}
\|\varrho_j\|_{W^{1, \infty}(\Rr^d\times (-1, 1))}\le C(1+N_s), \quad C=C(h, s, d),
\eq
whence  \eqref{est:Phij} follows. On the other hand, we have  $\Phi^{-1}_j(x, y)=(x, \ka_j(x, y))$ where 
\bq\label{def:kaj}
y= \varrho_j(x, z) \iff  z= \ka_j(x, y)\quad a.e.~ (x, z)\in \Rr^d\times (-1, 1).
\eq
Then  the relation $\ka_j(x, \varrho_j(x, z))=z$ yields
\bq\label{dkaj}
\p_y\ka_j(x, \varrho_j(x, z))=\frac{1}{\p_z\varrho_j(x, z))},\quad \p_x\ka_j(x, \varrho_j(x, z))=-\frac{\p_x\varrho_j(x, z)}{\p_z\varrho_j(x, z))}.
\eq
Thus,  in view of  \eqref{lower:varrho} and \eqref{est:varrho},  we obtain \eqref{est:Phi-1}. 
\end{proof}
\begin{lemm}\label{lemm:Psi}
Set 
\bq\label{def:Upsilon}
\Upsilon(x, y)=
\begin{cases}
\Phi_1\circ \Phi_2^{-1},\quad (x, y)\in \Omega_*,\\
(x, y),\quad (x, y)\in \Omega\setminus\Omega_*
\end{cases}
\eq
and 
\bq\label{def:M}
M=\frac{\na\Upsilon \na\Upsilon^t}{|\det \na \Upsilon|}.
\eq
Then,  $\Upsilon$ is a Lipschitz diffeomorphism on $\Omega$ and
\begin{align}\label{est:Jacobian}
&\frac{1}{C(1+N_s)}\le\det\na \Upsilon(x, y) \le C(1+N_s)\quad a.e.~ (x, y)\in \Omega,\\ \label{est:dUpsilon}
& \| \na\Upsilon\|_{L^\infty(\Omega)}+ \|  \na(\Upsilon^{-1}) \|_{L^\infty(\Omega)}\le  \cF(N_s).
\end{align}
Moreover, $M$ satisfies 
\begin{align}
\label{est:M}
&\| M-\mathrm{Id}\|_{L^\infty(\Omega)}\le \cF(N_s)\| \eta_\delta\|_{H^s},\\ \label{estM:H1}
&\| M-\mathrm{Id}\|_{L^2(\Omega)}\le \cF(N_s)\| \eta_\delta\|_{H^\mez}.
\end{align}
\end{lemm}
\begin{proof}
According to Lemma \ref{diffeos}, $\Upsilon$ is a Lipschitz diffeomorphism on $\Omega_*$. For $(x, y)\in \Omega\setminus\Omega_*$ we have $\Upsilon=\mathrm{Id}$, and hence $\Upsilon$ is a Lipschitz diffeomorphism on $\Omega$ and $M-\mathrm{Id}=0$. It thus suffices to consider $(x, y)\in \Omega_*$. On $\Omega_*$ we have $\Upsilon(x, y)=(x, \varrho_1(x, \ka_2(x, y))$, and so
\[
\na \Upsilon(x, y)=\begin{pmatrix}
1 & 0\\
a(x, y) & b(x, y)
\end{pmatrix}
\]
where
\[
a(x, y)=\p_x\varrho_1(x, \ka_2(x, y))+\p_z\varrho_1(x, \ka_2(x, y))\p_x\ka_2(x, y),\quad b(x, y)=\p_z\varrho_1(x, \ka_2(x, y))\p_y\ka_2(x, y).
\]
Using \eqref{dkaj} (with $j=2$) gives
\begin{align*}
&a(x, y)=\frac{\p_x\varrho_1(x, \ka_2(x, y))\p_z\varrho_2(x, \ka_2(x, y))-\p_x\varrho_2(x, \ka_2(x, y))\p_z\varrho_1(x, \ka_2(x, y))}{\p_z\varrho_2(x, \ka_2(x, y))},\\
& b(x, y)=\frac{\p_z\varrho_1(x, \ka_2(x, y))}{\p_z\varrho_2(x, \ka_2(x, y))}.
\end{align*}
In view of \eqref{lower:varrho} and \eqref{est:varrho} we obtain
\bq\label{est:det}
\frac{1}{C(1+N_s)}\le\det\na \Upsilon(x, y)=b(x, y) \le C(1+N_s)\quad a.e.~ (x, y)\in \Omega_*.
\eq
Next we compute 
\[
M-\mathrm{Id}=\frac{1}{b}\begin{pmatrix}
1-b& a\\
a& a^2+b(b-1)
\end{pmatrix}.
\]
Using the above formulas for $a$ and $b$ together with \eqref{diffeo} and \eqref{lower:varrho} we deduce that 
\begin{align*}
&\| (a, b-1)\|_{L^\infty(\Omega_*)}\le\cF(N_s)\| \eta_\delta\|_{H^s},\\
&\| (a, b-1)\|_{L^2(\Omega_*)}\le\cF(N_s)\| \eta_\delta\|_{H^\mez}.
\end{align*}
This combined with \eqref{est:det} leads to \eqref{est:dUpsilon}, \eqref{est:M} and \eqref{estM:H1}.
\end{proof}
\subsection{Proof of Proposition \ref{gamContLem}}
The proof proceeds in three steps. 

{\bf Step 1.} We first recall from  \eqref{JSystem} that $J^\pm v=f^\pm$ where $f^\pm$ solve 
\begin{equation}\label{JSys:1}
\begin{pmatrix}
\text{Id} & -\mathrm{Id}\\
\mu^+ G^-(\eta) & -\mu^- G^+(\eta) 
\end{pmatrix}
\begin{pmatrix}
f^- \\ f^+
\end{pmatrix} = \begin{pmatrix}
v \\ 0
\end{pmatrix}.
\end{equation}
From the definition of the Dirichlet-Neumann operator we see that  $f^\pm=q^\pm\vert_\Sigma$ where $q^\pm$ solve the two-phase elliptic problem 
\bq\label{elliptic:qpm}
\begin{cases}
\Delta q^\pm=0\quad\text{in}~\Omega^\pm,\\
q^--q^+=v \quad\text{on}~\Sigma,\\
\frac{\p_nq^-}{\mu^-}-\frac{\p_nq^+}{\mu^+}=0\quad\text{on}~\Sigma,\\
\p_{\nu^\pm}q^\pm=0\quad\text{on}~\Gamma^\pm.
\end{cases}
\eq
To remove the jump of $q$ at $\Sigma$ we take a function $\tt:\Omega\to \Rr$ satisfying  
\begin{align}\label{cd:tt}
&\tt(x, \eta(x))=-\mez v(x),\quad \tt\equiv 0\quad\text{near}~ \Gamma^\pm,\\ \label{est:tt}
& \| \tt\|_{\dot H^1(\Omega)}\le C(1+\| \eta\|_{W^{1, \infty}})\| v\|_{H^\mez(\Rr^d)},\quad C=C(d).
\end{align}
  Then, the solution of \eqref{elliptic:qpm} can be taken to be  $q^\pm:=(r\pm \tt)\vert_{\Omega^\pm}$ where $r\in \dot H^1(\Omega)$ solves 
 \bq\label{elliptic:rpm}
\begin{cases}
-\Delta r=\pm \Delta \tt\quad\text{in}~\Omega^\pm,\\
  \frac{\p_n r}{\mu^+}-\frac{\p_n r}{\mu^-}=-\p_n\tt(\frac{1}{\mu^+}+\frac{1}{\mu^-})\quad \text{on}~\Sigma,\\
\p_{\nu^\pm}r=0\quad\text{on}~\Gamma^\pm.
\end{cases}
\eq
A function  $\tt$ satisfying \eqref{cd:tt} and \eqref{est:tt} can be constructed as follows. Let $\varsigma(z):\Rr\to \Rr^+$  be a cutoff function that is identically $1$ for $|z|\le \mez$ and vanishes for $|z|\ge 1$. Set
\bq\label{def:tt}
\underline\tt(x, z)=-\mez\varsigma(z)e^{-|z|\langle D_x\rangle}v(x),\quad\theta(x, y) =\tt\big(x, \frac{y-\eta(x)}{h}\big).
\eq
Then,  $\tt(x, y)=0$ for $|y-\eta(x)|\ge h$, and hence $\tt_1\equiv 0$ near $\Gamma^\pm$ in view of the condition $\dist(\eta, \Gamma^\pm)>h$. Moreover, \eqref{est:tt} is satisfied. 

Integration by parts leads to the following variational form of \eqref{elliptic:rpm}:
\begin{equation}\label{weakform:2p}
\int_\Omega \Big(\frac{1_{\Omega^-}}{\mu^-}+ \frac{1_{\Omega^+}}{\mu^+}\Big)\nabla r \cdot \nabla \phi \ dxdy= \int_\Omega \Big(\frac{1_{\Omega^-}}{\mu^-}-\frac{1_{\Omega^+}}{\mu^+}\Big)\nabla \theta \cdot \nabla \phi \ dxdy,\quad\forall \phi\in \dot H^1(\Omega).
\end{equation}
For example, for $\varsigma(z):\Rr\to \Rr^+$  a cutoff function that is identically $1$ for $|z|\le \mez$ and vanishes for $|z|\ge 1$, we set
\bq\label{def:tt}
\underline\tt(x, z)=-\mez\varsigma(z)e^{-|z|\langle D_x\rangle}v(x),\quad\theta(x, y) =\tt\big(x, \frac{y-\eta(x)}{h}\big).
\eq
Then,  $\tt(x, y)=0$ for $|y-\eta(x)|\ge h$, and hence $\tt_1\equiv 0$ near $\Gamma^\pm$ in view of the condition $\dist(\eta, \Gamma^\pm)>h$. Moreover, we have 
\bq\label{est:tt}
\| \tt\|_{\dot H^1(\Omega)}\le C(1+\| \eta\|_{W^{1, \infty}})\| v\|_{H^\mez(\Rr^d)},\quad C=C(d).
\eq 
By virtue of the Lax-Milgram theorem, there exists a unique solution $r\in \dot H^1(\Omega)$ to \eqref{weakform:2p} which obeys the bound
\bq\label{est:r}
\| r\|_{\dot H^1(\Omega)}\le C(\mu^\pm)\| \tt\|_{\dot H^1(\Omega)}\le  C(1+\| \eta\|_{W^{1, \infty}})\| v\|_{H^\mez(\Rr^d)}.
\eq
Consequently, 
\bq\label{formula:Jpm}
J^\pm v=f^\pm=\mathrm{Tr}_{\Omega^\pm\to \Sigma}(r\pm \tt),
\eq
 and hence by the trace operation \eqref{Tr}, 
\bq
\| J^\pm v\|_{\wt H^\mez_\pm(\Rr^d)}\le \cF(N_s)\| v\|_{H^\mez(\Rr^d)}. 
\eq
 We note that $\tt$ defined by \eqref{def:tt} depends on $\eta$, and so does $r$. 

{\bf Step 2.} In this step we prove contraction estimates for $J^\pm_\delta= J^\pm_1-J^\pm_2$ in $\wt H^\mez_\pm$. Recall that $\Upsilon$ defined by \eqref{def:Upsilon} is a Lipschitz diffeomorphism on $\Omega$. Let $\tt_1$ be defined as in \eqref{def:tt} with $\eta=\eta_1$, and let $\tt_2=\tt_1\circ \Upsilon$. Let us check that $\tt_2$ obeys \eqref{cd:tt} and \eqref{est:tt} for $\eta=\eta_2$. Indeed, using the fact that $\Upsilon:\Sigma_2\to \Sigma_1$ we have 
\[
\tt_2(x, \eta_2(x))=\tt_1(x, \Upsilon(x, \eta_2(x)))= \tt_1(x, \eta_1(x))=-\mez v(x).
\]
On the other hand, since $\Upsilon\equiv \mathrm{Id}$ near $\Gamma^\pm$ and $\tt_1\equiv 0$ near $\Gamma^\pm$, we deduce that  $\tt_1\equiv 0$ near $\Gamma^\pm$. Finally, the bound \eqref{est:tt} for $\tt_2$ follows from \eqref{est:tt} for $\tt_1$ and the Lipschitz bound \eqref{est:dUpsilon} for $\Upsilon$.
 
 According to Step 1, we have
\[
J^\pm_j v=\mathrm{Tr}_{\Omega_j^\pm\to \Sigma_j}(r_j\pm \tt_j)
\]
where   $r_j\in \dot H^1(\Omega)$ satisfies 
\begin{equation}\label{weakform:2p}
\int_\Omega \Big(\frac{1_{\Omega^-}}{\mu^-}+ \frac{1_{\Omega^+}}{\mu^+}\Big)\nabla r_j \cdot \nabla \phi \ dxdy= \int_\Omega \Big(\frac{1_{\Omega^-}}{\mu^-}-\frac{1_{\Omega^+}}{\mu^+}\Big)\nabla \theta_j \cdot \nabla \phi \ dxdy\quad\forall \phi\in \dot H^1(\Omega).
\end{equation}
Set $\wt r_2  = r_2 \circ \Upsilon^{-1}$ and recall that $\theta_1=\theta_2\circ \Upsilon^{-1}$.  Combining \eqref{est:r} and \eqref{est:Jacobian} gives
\[
\| \wt r_2\|_{\dot H^1(\Omega)}\le \cF(N_s).
\]
Since map $\phi \mapsto \phi \circ \Upsilon^{-1}$ is an isomorphism on $\dot H^1 (\Omega)$, with $M = \nabla \Upsilon \nabla \Upsilon^t/|\det\nabla \Upsilon|$  we have for all $\phi \in \dot H^1(\Omega)$ that
\begin{equation}\label{weakform:2pnew}
\begin{aligned}
&\int_\Omega \Big(\frac{1_{\Omega_1^-}}{\mu^-}+ \frac{1_{\Omega_1^+}}{\mu^+}\Big) \nabla \tilde r_2 M \nabla \phi^t \ dxdy  = \int_\Omega \Big(\frac{1_{\Omega_1^-}}{\mu^-}-\frac{1_{\Omega_1^+}}{\mu^+}\Big)\nabla  \theta_1 M\nabla \phi ^t \ dxdy,
\end{aligned}
\end{equation}
where gradients of scalar functions are understood as row vectors, and the rows of the Jacobian matrix $\nabla \Upsilon$ are the gradients of each component of $\Upsilon$.  Taking the difference between \eqref{weakform:2pnew} with $j=2$ and \eqref{weakform:2p} with $j=1$, we obtain
\begin{equation}
\begin{aligned}
\int_\Omega \Big(\frac{1_{\Omega_1^-}}{\mu^-}+ \frac{1_{\Omega_1^+}}{\mu^+}\Big) \nabla (\wt r_2 -r_1)  \nabla \phi^t \ dxdy  &=  -\int_\Omega \Big(\frac{1_{\Omega_1^-}}{\mu^-}+\frac{1_{\Omega_1^+}}{\mu^+}\Big)\nabla  \wt r_2 (M-\mathrm{Id})\nabla \phi ^t \ dxdy\\
&\quad+\int_\Omega \Big(\frac{1_{\Omega_1^-}}{\mu^-}-\frac{1_{\Omega_1^+}}{\mu^+}\Big)\nabla  \theta_1 (M-\mathrm{Id})\nabla \phi ^t \ dxdy
\end{aligned}\label{diffVar}
\end{equation}
for all $\phi\in \dot H^1(\Omega)$. Setting $\phi=\wt r_2 -r_1$ and using the estimate \eqref{est:M} for $M-\mathrm{Id}$ in $L^\infty(\Omega)$ we obtain 
\bq
\label{est:rdelta:1}
\begin{aligned}
\| \wt r_2 -r_1\|_{\dot H^1(\Omega)}&\le C(\mu^\pm)\| M-\mathrm{Id}\|_{L^\infty(\Omega)}(\|\wt r_2\|_{\dot H^1(\Omega)}+\|\tt_1\|_{\dot H^1(\Omega)}) \\
&\le \cF(N_s)\| \eta_\delta\|_{H^s}\| v\|_{H^\mez}.
\end{aligned}
\eq
On the other hand, using \eqref{estM:H1} for $M-\mathrm{Id}$ in $L^2(\Omega)$ instead gives
\bq
\label{est:rdelta:2}
\begin{aligned}
\| \wt r_2 -r_1\|_{\dot H^1(\Omega)}&\le C(\mu^\pm)\| M-\mathrm{Id}\|_{L^2(\Omega)}(\|\na\wt r_2\|_{L^\infty(\Omega)}+\|\na\tt_1\|_{L^\infty(\Omega)}) \\
&\le \cF(N_s)\| \eta_\delta\|_{H^\mez}\| v\|_{H^s},
\end{aligned}
\eq
where in the last inequality we have used the fact that 
\[
\| \wt r_2\|_{L^\infty(\Omega)}+\|\tt_1\|_{L^\infty(\Omega)}\le \cF(N_s)\| v\|_{H^s}.
\]
Since $r_2=\wt r_2\circ\Upsilon$, $\tt_2= \tt\circ\Upsilon$, and $\Upsilon:\Omega_2^\pm\to \Omega_1^\pm,~\Sigma_2\to \Sigma_1$ we have 
\[
\mathrm{Tr}_{\Omega^\pm_2\to \Sigma_2}(r_2\pm \tt_2)=\mathrm{Tr}_{\Omega_1^\pm\to \Sigma_1}(\wt r_2\pm \tt_1),
\]
and hence
\begin{align*}
J^\pm_\delta v&=\mathrm{Tr}_{\Omega_1^\pm\to \Sigma_1}(r_1\pm \tt_1)-\mathrm{Tr}_{\Omega_2^\pm\to \Sigma_2}(r_2\pm \tt_2)\\
&=\mathrm{Tr}_{\Omega^\pm_1\to \Sigma_1}(r_1\pm \tt_1)-\mathrm{Tr}_{\Omega_1^\pm\to \Sigma_1}(\wt r_2\pm \tt_1)\\
&=\mathrm{Tr}_{\Omega_1^\pm\to \Sigma_1}(r_1-\wt r_2).
\end{align*}
In view of \eqref{est:rdelta:1} and \eqref{est:rdelta:2}, the trace operation \eqref{Tr} yields
\begin{align}\label{contra:J:low1}
&\| J^\pm_\delta v\|_{H^\mez_\pm}\le \cF(N_s)\| \eta_\delta\|_{H^s}\| v\|_{H^\mez},\\\label{contra:J:low2}
&\| J^\pm_\delta v\|_{H^\mez_\pm}\le \cF(N_s)\| \eta_\delta\|_{H^\mez}\| v\|_{H^s}.
\end{align}
For the proof of \eqref{gamCont12}, we shall only need \eqref{contra:J:low2}.

{\bf Step 3.}  We have $J^\pm_jv\equiv J^\pm(\eta_j)v=f^\pm_j$ where 
\bq\label{sys:fpmj}
\begin{cases}
f^-_j-f^+_j=v,\\
\frac{1}{\mu^-}G^-_jf_j^--\frac{1}{\mu^+}G^+_jf_j^+=0,
\end{cases}
\quad j=1, 2.
\eq
By taking differences we obtain $f^-_\delta=f^+_\delta$ and
\bq\label{eq:diffJpm}
\frac{1}{\mu^-} G^-_1  f^-_\delta-\frac{1}{\mu^+} G^+_1  f^+_\delta=\frac{1}{\mu^+}G^+_\delta  f^+_2-\frac{1}{\mu^-}G^-_\delta  f^- _2,
\eq
where we recall the notation $G^\pm_\delta=G^\pm_1-G^\pm_2$. Combining the contraction estimate \eqref{contraction:DN2}   with the continuity  \eqref{continuity:JGL} for $J^\pm$, 
we deduce that
\begin{align}
\label{est:Gdelta2}
&\|G^+_\delta  f^+_2\|_{H^{\sigma-1}} +\|G^-_\delta  f^-_2\|_{H^{\sigma-1}}\le \cF(N_s)\| \eta_\delta\|_{H^\sigma}\| v\|_{H^s}\quad\forall \sigma\in [\mez, s].
\end{align}
We take $\delta\in (0, \mez]$ satisfying $\sigma<s-1-\frac d2$. In light of the paralinearization \eqref{DNprincipal} for $G^\pm_1$, we have 
\bq\label{est:Jdelta:10}
\begin{aligned}
&\frac{1}{\mu^-} G^-_1  f^-_\delta -\frac{1}{\mu^+} G^+_1  f^+_\delta  = T_{\lambda_1}\big(\frac{1}{\mu^-} f^-_\delta  + \frac{1}{\mu^+} f^+_\delta \big) + R^1=\frac{\mu^++\mu^-}{\mu^-\mu^+}T_{\lambda_1} f^-_\delta+R^1,\\
&\| R^1\|_{H^{\nu-1+\delta}}\le \cF(N_s)\| f^-_\delta\|_{\wt H^\nu_-}\quad\forall \nu\in [\mez, s-\delta].
\end{aligned}
\eq
It follows from \eqref{eq:diffJpm},\eqref{est:Gdelta2} and \eqref{est:Jdelta:10} that if
\bq\label{cond:nu-sigma}
\nu\in [\mez, s-\delta],\quad \nu+\delta\le \sigma,\quad\sigma\in[\mez, s],
\eq
then
\begin{align}
&\| T_{\lambda_1} f^-_\delta\|_{H^{\nu-1+\delta}}\le\cF(N_s)\| f^-_\delta\|_{\wt H^\nu_-}+\cF(N_s)\| \eta_\delta\|_{H^\sigma}\| v\|_{H^s}. \label{bootstrap:Jdelta0}
\end{align}
Applying Lemma \ref{lemm:invertOp} we have
\[
\| \Psi(D) f^-_\delta\|_{H^{\nu+\delta}}\le\cF(\| \eta_1\|_{H^s})\| T_{\lambda_1} f^-_\delta\|_{H^{\nu+\delta-1}}+\cF(\| \eta_1\|_{H^s})\|f^-_\delta \|_{H^{1, \nu}}.
\]
But for $\tau>\mez$, 
\[
\| \cdot\|_{\wt H^{\tau}_\pm}\le C\| \Psi(D)\cdot\|_{H^\tau}+ C\|  \cdot\|_{\wt H^\mez_\pm},
\]
whence \eqref{bootstrap:Jdelta0} implies that
\begin{align}
&\|  f^-_\delta\|_{\wt H^{\nu+\delta}_-}\le\cF(N_s)\| f^-_\delta\|_{\wt H^\nu_-}+\cF(N_s)\| \eta_\delta\|_{H^\sigma}\| v\|_{H^s} \label{bootstrap:Jdelta}
\end{align}
provided that $\nu$ and $\sigma$ satisfy \eqref{cond:nu-sigma}. 

We note that \eqref{contra:J:low2} implies that
\bq\label{gamCont12:low}
\| f^-_\delta\|_{\wt H^\mez_-}\le \cF(N_s)\| \eta_\delta\|_{H^\mez}\| v\|_{H^s},
\eq
and hence  \eqref{gamCont12} holds for $\sigma=\mez$. Now fix $\sigma \in (\mez, s]$. We use \eqref{bootstrap:Jdelta}  to bootstrap the base estimate \eqref{gamCont12:low} to 
\bq\label{gamCont12:sigma}
\| f^-_\delta\|_{\wt H^\sigma_-}\le \cF(N_s)\| \eta_\delta\|_{H^\sigma}\| v\|_{H^s}.
\eq
Indeed, we can always  take $\delta$ smaller if necessary so that $\mez+\delta<\sigma$.  Plugging \eqref{gamCont12:low}  into \eqref{bootstrap:Jdelta} with $\nu=\mez$ yields \eqref{gamCont12:sigma} with $\mez+\delta$ in place of $\sigma$. Continuing this $n$ steps, $n$ being the greatest integer such that $\mez+n\delta\le \sigma$, we obtain  \eqref{gamCont12:sigma} for $\mez+n\delta$ in place of $\sigma$. This is justified since $\nu=\mez+(n-1)\delta$ satisfies \eqref{cond:nu-sigma}. Thus, for possibly one more step to gain $\sigma-(\mez+n\delta)$ derivative, we obtain \eqref{gamCont12:sigma}. The proof of \eqref{gamCont12} is complete. 

Finally, \eqref{gamCont11} can be proved similarly except that one uses the contraction estimate \eqref{contraction:DN3} to estimate $G^\pm_\delta  f^\pm_2$ in \eqref{est:Gdelta2}.


\section{Proof of Theorem \ref{theo:main}}\label{section:proof}
Let $s > 1+\frac{d}{2}$, $\mu^- > 0$, $\mu^+ \geq 0$, and $\mathfrak s \in (0, 1]$. Consider an initial datum $\eta_0\in H^s(\Rr^d)$ satisfying  
\begin{align}\label{RT:unilw}
&\inf_{x\in \Rr^d}   \mathrm{RT}(\eta_0)\ge 2\ma>0,\\\label{sep:initial}
&\dist(\eta_0, \Gamma^\pm)\ge 2h>0.
\end{align}
  Theorem \ref{theo:main} will be proved in Propositions \ref{prop:unitime}, \ref{prop:conv} and \ref{prop:conv2} below. Precisely, Proposition \ref{prop:unitime} establishes the uniform local well-posedness for the Muskat problem with surface tension belonging to any bounded set, say $\s\in (0, 1]$. Then, in Propositions \ref{prop:conv} and \ref{prop:conv2}, we prove that in appropriate topologies, $\eta^{(\s)}$ converges to $\eta$ with the rate $\sqrt\s$ and $\s$ respectively.
\begin{prop}\label{prop:unitime}
There exists a time $T_*>0$ depending only on $\| \eta_0\|_{H^s}$ and  $(\ma, h, s, \mu^\pm, \g)$  such that the following holds. For each $\s\in (0, 1]$, there exists a unique solution 
\bq\label{regularity:etas:0}
\eta^{(\s)}\in C([0, T_*]; H^s(\Rr^d))\cap L^2([0, T_*]; H^{s+\tdm}(\Rr^d))
\eq 
 to the Muskat problem with surface tension $\s$, $\eta^{(\s)}\vert_{t=0}=\eta_0$ and
 \bq\label{regularity:etas}
 \|\eta^{(\s)}\|_{L^\infty([0,T_*];H^s)}^2+ \|\eta^{(\s)}\|_{L^2([0,T_*];H^{s+\mez})}^2 + \mathfrak s \|\eta^{(\s)}\|_{L^2([0,T_*];H^{s+\tdm})}^2  \le  \mathcal F(\|\eta_0\|_{H^s}, \ma^{-1})
 \eq
for  some nondecreasing function $\cF:\Rr^+\times \Rr^+\to \Rr^+$  depending only on $(h, s, \mu^\pm, \g)$. Furthermore, for all $\s\in (0, 1]$ we have
\begin{align}\label{RT:etas}
&\inf_{t\in [0, T_*]}\inf_{x\in \Rr^d}   \mathrm{RT}(\eta^{(\s)}(t)) > \tdm\ma,\\ \label{seperation:etas}
&\inf_{t\in [0, T_*]}\dist(\eta^{(\s)}(t), \Gamma^\pm)> \tdm h.
\end{align}
\end{prop}
\begin{proof}
According to Theorems 1.2 and 1.3 in \cite{Ng1}, for each initial datum $\eta_0\in H^s$ satisfying $\dist(\eta_0, \Gamma^\pm)\ge 2h>0$ and for each $\s>0$, there exists $T_\s>0$ such that the Muskat problem has a unique solution
\[
\eta^{(\s)}\in C([0, T_\s]; H^s)\cap L^2([0, T_\s]; H^{s+\tdm})
\]
satisfying $\inf_{t\in T_\s}\dist(\eta^{(s)}(t),\Gamma^\pm)> \tdm h$. We stress the continuity in time of the $H^s$ norm of $\eta^{(\s)}$. Now we have in addition that $\eta_0$ satisfies the Rayleigh-Taylor condition \eqref{RT:unilw}. Thus, we define
\bq\label{def:T*s}
T^*_\s=\sup\{T\in (0, T_\s]: \inf_{t\in [0, T]}\inf_{x\in \Rr^d}   \mathrm{RT}(\eta^{(\s)}(t))> \tdm\ma\}.
\eq
 We shall prove that $T^*_\s>0$ for each $\s\in (0, 1]$ and  there exists $T_*>0$ such that $T^*_\s\ge T_*$ for all $s\in (0, 1]$.

{\bf Step 1.} We claim that there exist $\tt>0$ depending only on $s$, and $\cF_0: \Rr^+\to \Rr^+$ depending only on $(h, s, \mu^\pm, \g)$ such that
\bq\label{claim:A(t)}
\Big\| \lbb\mathfrak B(\eta^{(\s)}(t))J(\eta^{(\s)}(t))\rbb\eta^{(\s)}(t)- \lbb\mathfrak B(\eta_0)J(\eta_0)\rbb\eta_0\Big\|_{L^\infty(\Rr^d)}\le (t^{\frac{\tt}{2}}+ t^{\tt}) \cF_0(E_\s(t))
\eq
for all $t\le T_\s$, where
\bq\label{def:Es}
E_\s(t)=\| \eta^{(\s)}\|_{C([0,  t]; H^s)}^2+\s \| \eta^{(\s)}\|^2_{L^2([0,  t]; H^{s+\tdm})}.
\eq
Set 
\[
A(t)=\lbb\mathfrak B(\eta^{(\s)}(t))J(\eta^{(\s)}(t))\rbb\eta^{(\s)}(t)- \lbb\mathfrak B(\eta_0)J(\eta_0)\rbb \eta_0.
\]
The continuity properties \eqref{continuity:JGL} and \eqref{continuity:BV} of $J^\pm$ and $B^\pm$ imply that 
\bq\label{A(t):Hs}
\begin{aligned}
\| A(t)\|_{H^{s-1}}&\le \| \lbb\mathfrak B(\eta^{(\s)}(t))J(\eta^{(\s)}(t))\rbb \eta^{(\s)}(t)\|_{H^{s-1}}+\| \lbb\mathfrak B(\eta_0)J(\eta_0)\rbb\eta_0\|_{H^{s-1}}\\
&\le \cF\big(\| \eta^{(\s)}(t)\|_{H^s}\big)+\cF\big(\|\eta^{(\s)}(0)\|_{H^s}\big)\\
&\le \cF(E_\s(t)).
\end{aligned}
\eq
On the other hand, denoting $\mathfrak B^\pm(\eta^{(\s)}(t))=\frak B^\pm_t$ and  $J^\pm(\eta^{(\s)}(t))=J^\pm_t$,
 we can write 
\begin{align*}
A(t)&=\big(\mathfrak B^-_tJ^-_t-\frak B^+_tJ^+_t\big)\eta^{(\s)}(t)-\big(\mathfrak B^-_0J^-_0-\frak B^+_0J^+_0\big)\eta_0\\
&=(\frak B^-_t-\frak B^-_0)J^-_t\eta^{(\s)}(t)+\frak B^-_0(J^-_t-J^-_0)\eta^{(\s)}(t)+\frak B_0^-J_0^-(\eta^{(\s)}(t)-\eta_0)\\
&\qquad-(\frak B^+_t-\frak B^+_0)J^+_t\eta^{(\s)}(t)-\frak B^+_0(J^+_t-J^+_0)\eta^{(\s)}(t)-\frak B_0^+J_0^+(\eta^{(\s)}(t)-\eta_0).
\end{align*}
We treat the first two terms since the other terms are either similar or easier. The contraction estimate \eqref{contraction:DN2} with $\sigma=s-\mez$ gives
\[
\| G^\pm(\eta^{(\s)}(t))-G^\pm(\eta_0)\|_{\wt H^s_\pm \to H^{s-\tdm}}\le \cF(E_\s(t))\| \eta^{(\s)}(t)-\eta_0\|_{H^{s-\mez}}.
\]
From this and the definition of $\frak B^\pm$ it is easy to prove that 
\[
\| \frak B^+_t-\frak B^+_0\|_{ \wt H^s_\pm \to H^{s-\tdm}}\le \cF(E_\s(t))\| \eta^{(\s)}(t)-\eta_0\|_{H^{s-\mez}}.
\]
Then recalling the continuity of $J^\pm$ from $H^{s-\mez}\to \wt H^{s-\mez}_\pm$ we obtain 
\[
\| (\frak B^-_t-\frak B^-_0)J^-_t\eta^{(\s)}(t)\|_{H^{s-\tdm}}\le \cF(E_\s(t))\| \eta^{(\s)}(t)-\eta_0\|_{H^{s-\mez}}.
\]
Next regarding $\frak B^-_0(J^-_t-J^-_0)\eta^{(\s)}(t)$ we use the contraction estimate \eqref{gamCont12} with $\sigma=s-\mez$
\[
\| J ^-_t-J^-_0 \|_{ H^s\to \wt H^{s-\mez}_\pm} \leq \mathcal F(E_\s(t))\|\eta^{(\s)}(t)-\eta_0\|_{H^{s-\mez}}
\]
together with the continuity  \eqref{continuity:BV} to have 
\[
\| \frak B^-_0(J^-_t-J^-_0)\eta^{(\s)}(t)\|_{H^{s-\tdm}} \leq \mathcal F(E_\s(t))\|\eta^{(\s)}(t)-\eta_0\|_{H^{s-\mez}}.
\]
Therefore, we arrive at 
\bq\label{A(t):Hs-tdm}
\| A(t)\|_{H^{s-\tdm}} \leq \mathcal F(E_\s(t))\|\eta^{(\s)}(t)-\eta_0\|_{H^{s-\mez}}.
\eq
To bound $\|\eta^{(\s)}(t)-\eta_0\|_{H^{s-\mez}}$ we first use  the mean-value theorem and equation \eqref{evol} to have
\begin{align*}
\| \eta^{(\s)}(t)-\eta_0\|_{H^{-\mez}}&\le \int_0^t \| \p_t\eta^{(\s)}(\tau)\|_{H^{-\mez}}d\tau\\
& \le \cF(E_\s(t))\int_0^t \g\| \eta^{(\s)}(\tau)\|_{H^\mez}+\s\| H\big(\eta^{(\s)}(\tau)\big)\|_{H^\mez}d\tau \\
&\le  \cF(E_\s(t)) \Big(t\g\| \eta^{(\s)}\|_{C([0, t]; H^s)}+t^\mez\s\| \eta^{(\s)}\|_{L^2([0, t]; H^{s+\tdm})}\Big)\\
&\le  (t^\mez+t)\cF(E_\s(t)) \Big(\g\| \eta^{(\s)}\|_{C([0, t]; H^s)}+\sqrt\s\| \eta^{(\s)}\|_{L^2([0, t]; H^{s+\tdm})}\Big)
\end{align*}
for all $\s\in (0, 1]$ and $t\le T_\s$. Interpolating this with the obvious bound $\| \eta^{(\s)}(t)-\eta_0\|_{H^s}\le \cF(M_s(t))$ gives 
\bq\label{etat-eta0}
\|\eta^{(\s)}(t)-\eta_0\|_{H^{s-\mez}}\le  (t^{\frac{\tt_0}{2}}+ t^{\tt_0})\cF(E_\s(t)) 
\eq
for some $\tt_0\in (0, 1)$ and $\cF$ depending only on $(h, s, \mu^\pm, \g)$. Then in view of \eqref{A(t):Hs-tdm}, this implies
\[
\| A(t)\|_{H^{s-\tdm}}\le  (t^{\frac{\tt_0}{2}}+ t^{\tt_0})\cF(E_\s(t)).
\]
Fixing $s'\in (\max\{1+\frac d2, s-\tdm\}, s)$ and interpolating this with the $H^s$ bound \eqref{A(t):Hs} we obtain 
\[
\| A(t)\|_{H^{s'-1}}\le(t^{\frac{\tt}{2}}+ t^{\tt})\cF(E_\s(t))
\]
for some $\tt\in (0, 1)$. Using the embedding $H^{s'-1}\subset L^\infty(\Rr^d)$ we conclude the proof of \eqref{claim:A(t)}.

{\bf Step 2.} We note that \eqref{claim:A(t)} implies the continuity of 
\[
[0, T_\s]\ni t\mapsto \inf_{x\in \Rr^d}   \mathrm{RT}(\eta^{(\s)}(t)).
\]
 Thus, in view of the definition \eqref{def:T*s} and the initial condition \eqref{RT:unilw}, we have  $T^*_\s>0$ for all $\s \in (0, 1]$.

By the definition of $T^*_\s$, conditions \eqref{Regcond}, \eqref{RTcond} and \eqref{Sepcond} in Proposition \ref{apriori} are satisfied for all $T\le T^*_\s$.  Thus, the estimates \eqref{uni:Hs} and \eqref{uni:dissipation}  imply   the existence of a strictly increasing $\cF_2:\Rr^+\times\Rr^+\to \Rr^+$ depending only on $(h, s, \mu^\pm, \g)$ such that
\bq
\cF_2(m, 0)>m^2\quad\forall m>0
\eq
and
\bq\label{apriori:Es}
 M_\s(T)\le \cF_2\Big( \|\eta_0\|_{H^s} +T\mathcal F_2\big(M_\s(T),\mathfrak a^{-1}\big),\mathfrak a^{-1}\Big)
\eq
for all $\s\in (0, 1]$ and $T\le T_\s^*$, where 
\[
M_\s(T)
=\| \eta^{(\s)}\|_{C([0,  T]; H^s)}^2+\| \eta^{(\s)}\|_{L^2([0, T]; H^{s+\mez})}^2+\s \| \eta^{(\s)}\|^2_{L^2([0,  T]; H^{s+\tdm})}.
\]
  Set 
  \bq
 T_2= \frac{\|\eta_0\|_{H^s}}{2\mathcal F_2\Big(\cF_2\big( 2\|\eta_0\|_{H^s},\mathfrak a^{-1}\big),\mathfrak a^{-1}\Big)}
  \eq
independent of $\s$. We claim that 
\bq\label{Es-T2}
M_\s(T)\le K_0:=\cF_2\big(2\|\eta_0\|_{H^s}, \ma^{-1}\big)\quad\forall T\le \min\{T_2, T^*_\s\},\quad\forall \s\in (0, 1].
\eq
Assume not, then there exists $\s_0\in (0, 1]$ and $T_3\le \min\{T_2, T^*_{\s_0}\}$ such that $M_{\s_0}(T_3)> K_0$. Since $M_{\s_0}(0)=\|\eta_0\|_{H^s}^2<K_0$, the continuity of $T\mapsto E_{\s_0}(T)$ then yields the existence of $T_4\in (0, T_3)$ such that $M_{\s_0}(T_4)= K_0$. Consequently, at $T=T_4$, \eqref{apriori:Es} gives
\begin{align*}
  \cF_2\big( 2\|\eta_0\|_{H^s},\mathfrak a^{-1}\big)&\le \cF_2\Big( \|\eta_0\|_{H^s} +T_4\mathcal F_2\Big(\cF_2\big( 2\|\eta_0\|_{H^s},\mathfrak a^{-1}\big),\mathfrak a^{-1}\Big),\mathfrak a^{-1}\Big)\\
  &\le \cF_2\Big( \|\eta_0\|_{H^s} +T_2\mathcal F_2\Big(\cF_2\big( 2\|\eta_0\|_{H^s},\mathfrak a^{-1}\big),\mathfrak a^{-1}\Big),\mathfrak a^{-1}\Big)\\
    &\le \cF_2\Big( \tdm\|\eta_0\|_{H^s},\mathfrak a^{-1}\Big),
  \end{align*}
where we have used the definition of $T_2$ in the last inequality. This contradicts the fact that $\cF_2$ was chosen to be strictly increasing. 

Now for all $T\le  \min\{T_2, T^*_\s\}$, we use \eqref{claim:A(t)} and the fact that $E_\s(\cdot)\le M_\s(\cdot)$ to have 
\begin{align*}
\inf_{t\in [0, T]}\inf_{x\in \Rr^d} \mathrm{RT}(\eta^{(\s)}(t))&\ge \inf_{x\in \Rr^d} \mathrm{RT}(\eta_0)- (T^{\frac{\tt}{2}}+ T^{\tt}) \cF_0(M_\s(T))\\
&\ge 2\ma - (T^{\frac{\tt}{2}}+ T^{\tt}) \cF_0(K_0).
\end{align*}
Choosing $T_*\le T_2$ sufficiently small so that 
\bq\label{chooseT*:2}
(T_*^{\frac{\tt}{2}}+ T_*^{\tt}) \cF_0(K_0)<\mez\ma
\eq
we obtain
\bq\label{RT:T1}
\inf_{t\in [0, T]}\inf_{x\in \Rr^d} \mathrm{RT}(\eta^{(\s)}(t))>\tdm\ma \quad\forall T\le \min\{T_*, T_\s^*\}.
\eq
Clearly, $T_*$ is independent of $\s$. Moreover, \eqref{RT:T1} and the definition of $T_\s^*$ show that $T_*\le T^*_\s$ for all $\s\in (0, 1]$. Finally, since $T_*\le \min\{T_2, T_\s^*\}$, \eqref{Es-T2} and the  definition of $T^*_\s$  guarantee that the estimates \eqref{regularity:etas}, \eqref{RT:etas} and \eqref{seperation:etas} hold true. 
\end{proof}

\begin{prop}\label{prop:conv}
There exists $\cF:\Rr^+\times \Rr^+\to \Rr^+$ depending only on $(h, s, \mu^\pm, \g)$ such that
\bq\label{proof:finalest}
\begin{aligned}
&\|\eta^{(\s_1)}-\eta^{(\s_2)}\|_{L^\infty([0,T_*];H^{s-1})}^2+ \|\eta^{(\s_1)}-\eta^{(\s_2)}\|_{L^2([0,T_*];H^{s-\mez})}^2 \le   (\s_1+\s_2)\mathcal F\big(\| \eta_0\|_{H^s}, \ma^{-1}\big)
\end{aligned}
\eq
for all $\s_1$ and $\s_2$ in $(0, 1]$.
\end{prop}
\begin{proof}
Denote $\eta_j=\eta^{(\s_j)}$, $j=1, 2$, and $\eta_\delta=\eta_1-\eta_2$ which exists on $[0, T_*]$. We fix  $\delta\in \big(0, \min(\mez, s-1-\frac d2)\big)$. From \eqref{evol} we have that $\eta_\delta$ evolves according to
\begin{equation} \label{deltSys}
(\mu^+ + \mu^-) \partial_t \eta_\delta = -\mathfrak g (L_\delta \eta_1+ L_2\eta_\delta)- \mathfrak s_1 L_1 H(\eta_1)+\mathfrak s_2 L_2 H(\eta_2).
\end{equation}
 By \eqref{continuity:JGL} and \eqref{boundH}, 
\bq\label{L1H1}
\| L_j H(\eta_j)\|_{H^{s-\tdm}}\le \cF(\|\eta_j\|_{H^s})\| H(\eta_j)\|_{H^{s-\mez}}\le  \cF(\|\eta_j\|_{H^s})\| \eta_j\|_{H^{s+\tdm}}.
\eq
 We now paralinearize $L_2$ and $L_\delta$.  Applying \eqref{Lprincipal} with $\sigma=s-\mez-\delta$ gives
\bq\label{lin:L2}
L_2 \eta_\delta = T_{\lambda_2}\eta_\delta+  O_{H^{s-\mez-\delta} \to H^{s -\tdm}}(\mathcal F(N_s))\eta_\delta.
\eq
As for $L_\delta$ we write 
\begin{align}\label{expand:Ldelta}
L_\delta &= \sum_{\pm = +,-} G_\delta^\pm  J _1^\pm + G^\pm_2  J _\delta^\pm.
\end{align}
Using \eqref{DNprincipal} at $\sigma = s-\mez-\delta$, we have
\begin{align*}
&\sum_{\pm = +,-} G_2^\pm  J _\delta^\pm = T_{\lambda_2}(J _\delta^--J^+_\delta)   + \sum_{\pm = +,-} O_{\wt H^{s-\mez-\delta}_\pm\to H^{s-\tdm}}(\mathcal F(N_s)) J _\delta^\pm.
\end{align*}
Recall that $N_s$ is given by \eqref{Ns}. However, 
\[
J _\delta^--J^+_\delta=(J_1^--J_1^+)-(J_2^--J_2^+)= \text{Id} - \text{Id} = 0
\]
 and by virtue of Proposition \ref{gamContLem} (with $\sigma=s-\mez-\delta)$, 
 \begin{align*}
 &J^\pm_\delta=O_{H^s\to \wt H^{s-\mez-\delta}_\pm}(\| \eta_\delta\|_{H^{s-\mez-\delta}}).
 \end{align*}
 We thus obtain
\begin{align}\label{proof:Gdelta1}
&\sum_{\pm = +,-}G_2^\pm  J _\delta^\pm = O_{H^{s} \to H^{s -\tdm}}\big(\mathcal F(N_s)\|\eta_\delta\|_{H^{s-\mez-\delta}}\big).
\end{align}
As for $G_\delta^\pm  J _1^\pm$, we apply  Theorem \ref{DNcont} with $\sigma=s-\mez$  and \eqref{continuity:JGL} to have
\begin{align*}
\mp G^\pm_\delta  J _1^\pm f & = -T_{\lambda_1 \mathfrak B^\pm_1  J _1^\pm f}\eta_\delta  -T_{\mathfrak V^\pm_1  J _1^\pm f} \cdot\nabla  \eta_\delta + O_{H^{s} \to H^{s -\tdm}}\big(\mathcal F(N_s)\|\eta_\delta\|_{H^{s-\mez - \delta}}\big)f,
\end{align*}
and hence
\[
\sum_{\pm = +,-} G_\delta^\pm  J _1^\pm f=-T_{ \lambda_1 \lbb \mathfrak B_1  J _1\rbb f}\eta_\delta  -T_{\lbb \mathfrak V_1  J _1\rbb f} \cdot\nabla  \eta_\delta + O_{H^s \to H^{s-\tdm}}\big(\mathcal F(N_s) \|\eta_\delta\|_{H^{s-\mez-\delta}}\big)f.
\]
Combining this with  \eqref{proof:Gdelta1} yields
\bq\label{Ldelt}
L_\delta f =- T_{ \lambda_1 \lbb \mathfrak B_1  J _1\rbb f}\eta_\delta  -T_{\lbb \mathfrak V_1  J _1\rbb f} \cdot\nabla  \eta_\delta + O_{H^s \to H^{s-\tdm}}\big(\mathcal F(N_s) \|\eta_\delta\|_{H^{s-\mez-\delta}}\big)f.
\eq
From  \eqref{lin:L2} and \eqref{Ldelt} we have 
\begin{align*}
L_2\eta_\delta+L_\delta \eta_1&=T_{\lambda_2}\eta_\delta- T_{ \lambda_1 \lbb \mathfrak B_1  J _1\rbb \eta_1}\eta_\delta  -T_{\lbb \mathfrak V_1  J _1\rbb \eta_1} \cdot\nabla  \eta_\delta \\
&\qquad + O_{H^{s-\mez-\delta} \to H^{s -\tdm}}(\mathcal F(N_s))\eta_\delta+ O_{H^s \to H^{s-\tdm}}\big(\mathcal F(N_s) \|\eta_\delta\|_{H^{s-\mez-\delta}}\big)\eta_1
\end{align*}
Interchanging $\eta_1$ and $\eta_2$ gives 
\begin{align*}
L_1\eta_\delta+L_\delta \eta_2&=T_{\lambda_1}\eta_\delta- T_{ \lambda_2 \lbb \mathfrak B_2  J _2\rbb \eta_2}\eta_\delta  -T_{\lbb \mathfrak V_2  J _2\rbb \eta_2} \cdot\nabla  \eta_\delta \\
&\qquad + O_{H^{s-\mez-\delta} \to H^{s -\tdm}}(\mathcal F(N_s))\eta_\delta+O_{H^s \to H^{s-\tdm}}\big(\mathcal F(N_s) \|\eta_\delta\|_{H^{s-\mez-\delta}}\big)\eta_2.
\end{align*}
But $L_1\eta_\delta+L_\delta \eta_2=L_2\eta_\delta+L_\delta \eta_1$, thus taking the average of the above identities yields
\begin{align*}
L_2\eta_\delta+L_\delta \eta_1&= T_{ (\lambda (1-\lbb \mathfrak B  J \rbb \eta))_\alpha}\eta_\delta  -T_{(\lbb \mathfrak V  J\rbb \eta)_\alpha} \cdot\nabla  \eta_\delta \\
&\qquad + O_{H^{s-\mez-\delta} \to H^{s -\tdm}}(\mathcal F(N_s))\eta_\delta+ O_{H^s \to H^{s-\tdm}}\big(\mathcal F(N_s) \|\eta_\delta\|_{H^{s-\mez-\delta}}\big)(\eta_1+\eta_2),
\end{align*}
where 
\[
\big(\lambda (1-\lbb \mathfrak B  J \rbb \eta)\big)_\alpha=\mez\Big(\lambda_1 (1-\lbb \mathfrak B_1  J_1 \rbb \eta_1)+\lambda_2 (1-\lbb \mathfrak B_2  J_2 \rbb \eta_2)\Big)
\]
and similarly for $(\lbb \mathfrak V  J\rbb \eta)_\alpha$. It then follows from \eqref{deltSys} and \eqref{L1H1} that
\bq\label{proof:lin1}
\begin{aligned}
&(\mu^+ + \mu^-) \partial_t \eta_\delta =- \g T_{ (\lambda (1-\lbb \mathfrak B  J \rbb \eta))_\alpha}\eta_\delta  +\g T_{(\lbb \mathfrak V  J\rbb \eta)_\alpha} \cdot\nabla  \eta_\delta+\cR_1,\\
&\| \cR_1\|_{H^{s-\tdm}}\le \g\cF(N_s)\|\eta_\delta\|_{H^{s-\mez-\delta}}+(\s_1+\s_2)\cF(N_s)N_{s+\tdm}.
\end{aligned}
\eq
 Next we perform $H^{s-1}$ energy estimate for \eqref{proof:lin1}. Introduce $\eta_{\delta, s-1}=\langle D\rangle^{s-1}\eta_\delta$. Upon commuting \eqref{proof:lin1} with $\langle D\rangle^{s-1}$ and applying Theorem \ref{theo:sc}  we arrive at 
 \bq\label{proof:lin3}
\begin{aligned}
&(\mu^+ + \mu^-) \partial_t \eta_{\delta, s-1} =-\g T_{ (\lambda (1-\lbb \mathfrak B  J \rbb \eta))_\alpha}\eta_{\delta, s-1} +\g i\mathrm{Re}\big(T_{(\lbb \mathfrak V  J\rbb \eta)_\alpha\cdot \xi}\big) \eta_{\delta, s-1}+\cR,\\
&\| \cR\|_{H^{-\mez}}\le\g\cF(N_s)\|\eta_{\delta, s-1}\|_{H^{\mez-\delta}}+\cF(N_s)\sum_{j=1}^2\s_j\| \eta_j\|_{H^{s+\tdm}},
\end{aligned}
\eq
where $\cF$ depends only on $(h, s, \mu^\pm)$. Moreover, the uniform estimate \eqref {regularity:etas} implies that 
\bq\label{est:Ns}
N_s\le \cF(\| \eta_0\|_{H^s}, \ma^{-1}).
\eq
Testing  \eqref{proof:lin3} against $\eta_{\delta, s-1}$ yields
 \begin{align*}
\frac{(\mu^+ + \mu^-)}{2}\frac{d}{dt}\|\eta_{\delta, s-1} \|_{L^2}^2  =& - \mathfrak g \big( T_{  (\lambda (1-\lbb \mathfrak B  J \rbb \eta))_\alpha}\eta_{\delta, s-1},\eta_{\delta, s-1}\big)_{L^2}
+(\cR, \eta_{\delta, s-1})_{L^2},
\end{align*}
where we have used the fact that $i\mathrm{Re}\big(T_{(\lbb \mathfrak V  J\rbb \eta)_\alpha\cdot \xi}\big)$ is skew-adjoint.

From  \eqref{RT:etas}  we have that $ (\lambda (1-\lbb \mathfrak B  J \rbb \eta))_\alpha$ is an elliptic symbol in $\Gamma^1_\delta$:
\[
  (\lambda (1-\lbb \mathfrak B  J \rbb \eta))_\alpha\ge \ma |\xi|,\quad  M^1_\delta\big( (\lambda (1-\lbb \mathfrak B  J \rbb \eta))_\alpha\big)\le\cF(N_s, \ma^{-1}),
\]
where $\cF$ depends only on $(h, s, \mu^\pm, \g)$. Applying the  G\aa rding inequality \eqref{wgi}, we have
\[
\| \eta_{\delta, s-1} \|^2_{H^\mez}\le \cF(\| \eta_0\|_{H^s}, \ma^{-1})\Big(\big( T_{  (\lambda (1-\lbb \mathfrak B  J \rbb \eta))_\alpha}\eta_{\delta, s-1},\eta_{\delta, s-1}\big)_{L^2}+\| \eta_{\delta, s-1}\|_{H^\mez}\|  \eta_{\delta, s-1}\|_{H^{\mez-\delta}} \Big).
\]
This combined with the estimate for $\cR$ in \eqref{proof:lin3} and \eqref{est:Ns} implies
\bq\label{eest:etadelta}
\begin{aligned}
&\frac{(\mu^+ + \mu^-)}{2}\frac{d}{dt}\|\eta_{\delta, s-1} \|_{L^2}^2  +\frac{\g}{\mathcal{Q}_{T_*}}\| \eta_{\delta, s-1} \|^2_{H^\mez}\\
&\quad\le \g\mathcal{Q}_{T_*}\|\eta_{\delta, s-1}\|_{H^{\mez-\delta}}\| \eta_{\delta, s-1}\|_{H^{\mez}}+\mathcal{Q}_{T_*}\| \eta_{\delta, s-1}\|_{H^{\mez}}\sum_{j=1}^2\s_j\| \eta_j\|_{H^{s+\tdm}},
\end{aligned}
\eq
where 
\bq\label{Q:conv}
\mathcal{Q}_{T_*}=\cF\big(\| \eta_0\|_{H^s}, \ma^{-1}\big),
\eq
$\cF$ depending only on $(h, s, \mu^\pm, \g)$. Using interpolation and Young's inequality, we have
\bq\label{interpolate:conv}
\begin{aligned}
&\|\eta_{\delta, s-1}\|_{H^{\mez-\delta}}\|\eta_s\|_{H^{\mez}} \leq \|\eta_{\delta, s-1}\|_{L^2}^{2\delta}\|\eta_{\delta, s-1}\|_{H^{\mez}}^{2(1-\delta)} \leq (10\mathcal{Q}_{T_*})^{\frac{2(1-\delta)}{\delta}}\|\eta_{\delta, s-1}\|_{L^2}^2 + \frac{1}{100\mathcal Q^2_{T_*}}\|\eta_{\delta, s-1}\|_{H^{\mez}}^2,\\
& \s_j\| \eta_j\|_{H^{s+\tdm}}\| \eta_{\delta, s-1}\|_{H^{\mez}}\le \s_j^2100\mathcal{Q}^2_{T_*}\| \eta_j\|_{H^{s+\tdm}}^2+ \frac{1}{100\mathcal Q^2_{T_*}}\|\eta_{\delta, s-1}\|_{H^{\mez}}^2.
\end{aligned}
\eq
Thus, for possibly a larger $\cF$ in  $\mathcal{Q}_{T_*}$, we obtain 
\bq\label{diff:ineq1}
\begin{aligned}
&\frac{(\mu^+ + \mu^-)}{2}\frac{d}{dt}\|\eta_{\delta, s-1} \|_{L^2}^2  +\frac{\g}{\mathcal{Q}_{T_*}}\| \eta_{\delta, s-1}\|_{H^{\mez}}^2\le \g\mathcal{Q}_{T_1}\| \eta_{\delta, s-1}\|^2_{L^2}+\mathcal{Q}_{T_*} \sum_{j=1}^2\s_j^2\| \eta_j\|_{H^{s+\tdm}}^2.
\end{aligned}
\eq
Finally, since $\eta_\delta\vert_{t=0}=0$ and by \eqref{regularity:etas}
\[
\mathfrak s_j \int_0^{T_*}\|\eta_j\|^2_{H^{s+\tdm}}dt  \le  \mathcal F(\|\eta_0\|_{H^s}, \ma^{-1}), 
\]
 an application of Gr\"onwall's lemma leads to the  estimate  \eqref{proof:finalest}.
\end{proof}
Now let $\s_n\to 0$ and rename $\eta_n=\eta^{(\s_n)}$ solution to the Muskat problem with surface tension $\s_n$ on $[0, T_*]$. The uniform estimates in \eqref{regularity:etas} show that along a subsequence 
$\eta_n$ converges weakly-* to 
\bq\label{regeta:Hs}
\eta\in L^\infty([0, T_*]; H^s)\cap L^2([0, T_*]; H^{s+\mez})
\eq
together with the bounds
 \begin{align}\label{bound:limit1}
 &\|\eta\|_{L^\infty([0,T_*];H^s)}+ \|\eta\|_{L^2([0,T_*]; H^{s+\mez})}  \le  \mathcal F(\|\eta_0\|_{H^s}, \ma^{-1}).
  \end{align}
The estimate \eqref{proof:finalest} implies that $(\eta_n)_n$ is a Cauchy sequence in $C([0, T_*]; H^{s-1})\cap L^2([0, T_*]; H^{s-\mez})$. Therefore, 
\bq\label{converge:Z}
\eta_n\to \eta \quad\text{in}~C([0, T_*]; H^{s-1})\cap L^2([0, T_*]; H^{s-\mez});
\eq
in particular, $\eta\vert_{t=0}=\eta_0$. Moreover, by interpolating between $L^\infty_t H^s_x$ and $C_t H^{s-1}$, we deduce that $\eta\in C([0, T_*]; H^{s'})$ for all $s'<s$.   Since $\eta_n\to \eta$ in $C_t H^{s-1}_x\subset C_tL^\infty_x$,  \eqref{seperation:etas} gives
 \bq \label{bound:limit2}
\inf_{t\in [0, T_*]}\dist(\eta(t), \Gamma^\pm)\ge \tdm h.
\eq
\begin{lemm}
$\eta$ is a solution on $[0, T_*]$ of the Muskat problem without surface tension with initial data $\eta_0$.  
\end{lemm}
\begin{proof}
 For each $n$, we have from \eqref{evol} that
\begin{align}
 \partial_t\eta_n +\frac{1}{\mu^+ + \mu^-}L(\eta_n)( \mathfrak g \eta_n  + \mathfrak s H(\eta_n))= 0
\end{align}
For any compactly supported test function $\varphi\in C^\infty((0, T_*)\times \Rr^d)$, we have
\bq
\int_0^{T_*}\eta_n\p_t\varphi dxdt=\frac{1}{\mu^+ + \mu^-}\int_0^{T_*}\int_{\Rr^d}\varphi L(\eta_n)( \mathfrak g \eta_n  + \mathfrak s H(\eta_n))dxdt.
\eq
Clearly, \eqref{converge:Z} implies that 
\[
\int_0^{T_*}\eta_n\p_t\varphi dxdt\to \int_0^{T_*}\eta\p_t\varphi dxdt.
\] 
The continuity \eqref{continuity:JGL} of $L$ combined with \eqref{boundH} and the uniform bound \eqref{regularity:etas} yields 
\begin{align*}
\s_n \| L(\eta_n)H(\eta_n)\|_{L^2([0, T_*];  H^{s-\tdm})}&\le \cF(\|\eta_n\|_{L^\infty([0, T_*]; H^s)})\s_n\| \eta_n\|_{L^2([0, T_*]; H^{s+\tdm})}\\
&\les \sqrt{\s_n}\cF(\|\eta_0\|_{H^s}, \ma^{-1}).
\end{align*}
Since $s-\tdm>0$, this implies 
\begin{align*}
\left|\int_0^{T_*}\int_{\Rr^d}\varphi L(\eta_n) \big(\s_nH(\eta_n)\big)dxdt\right|&\le \| \varphi\|_{L^2_{x, t}}\| L(\eta_n) \big(\s_nH(\eta_n)\big)\|_{L^2_{x, t}}\\
&\les  \sqrt{\s_n}\| \varphi\|_{L^2_{x, t}}\cF(\|\eta_0\|_{H^s}, \ma^{-1})\to 0.
\end{align*}
Next we write 
\begin{align*}
L(\eta_n)\eta_n-L(\eta)\eta&=\big(L(\eta_n)-L(\eta)\big)\eta_n+L(\eta_n)(\eta_n-\eta)\\
&=\frac{\mu^++\mu^-}{\mu^-}\big(G^-(\eta_n)J^-(\eta_n)-G^-(\eta)J^-(\eta)\big)\eta_n+L(\eta_n)(\eta_n-\eta)\\
&=\frac{\mu^++\mu^-}{\mu^-}\Big\{\big(G^-(\eta_n)-G^-(\eta)\big)J^-(\eta_n)\eta_n-G^-(\eta)\big(J^-(\eta_n)-J^-(\eta)\big)\eta_n\Big\}\\
&\qquad+L(\eta_n)(\eta_n-\eta).
\end{align*}
Combining  \eqref{contraction:DN2} and \eqref{continuity:JGL} we obtain
\[
\|\big(G^-(\eta_n)-G^-(\eta)\big)J^-(\eta_n)\eta_n\|_{L^2([0, T_*];  H^{s-\tdm})}\le \cF(\|(\eta_n, \eta)\|_{L^\infty([0, T_*]; H^s)})\| \eta_n-\eta\|_{L^2([0, T_*]; H^{s-\mez})}.
\] 
On the other hand, \eqref{continuity:JGL} and \eqref{gamCont12} yield
\[
\| G^-(\eta)\big(J^-(\eta_n)-J^-(\eta)\big)\eta_n\|_{L^2([0, T_*]; H^{s-\tdm})}\le \cF(\|(\eta_n, \eta)\|_{L^\infty([0, T_*]; H^s)})\| \eta_n-\eta\|_{L^2([0, T_*]; H^{s-\mez})}.
\]
Finally, by \eqref{continuity:JGL} we have
\[
\| L(\eta_n)(\eta_n-\eta)\|_{L^2([0, T_*]; H^{s-\tdm})}\le \cF(\|\eta_n\|_{L^\infty([0, T_*]; H^s)})\| \eta_n-\eta\|_{L^2([0, T_*]; H^{s-\mez})}.
\]
Putting together the above considerations, we obtain
\[
\left| \int_0^{T_*}\int_{\Rr^d}\varphi\big(L(\eta_n)\eta_n-L(\eta)\eta\big)\right|\le \| \varphi\|_{L^2_{x,t}}\cF(\|(\eta_n, \eta)\|_{L^\infty([0, T_*]; H^s)})\| \eta_n-\eta\|_{L^2([0, T_*]; H^{s-\mez})}\to 0
\]
by virtue of the strong convergence \eqref{converge:Z} and the uniform $H^s$ bound in  \eqref{regularity:etas}. We have proved that 
\[
\int_0^{T_*}\eta\p_t\varphi dxdt=\frac{1}{\mu^+ + \mu^-}\int_0^{T_*}\int_{\Rr^d}\varphi L(\eta)(\mathfrak g \eta )dxdt
\]
for all compactly supported smooth test functions $\varphi$. Therefore, $\eta$ is a solution on $[0, T_*]$ of the Muskat problem without surface tension.
 \end{proof}
 \begin{lemm} We have
 \bq\label{bound:limit3}
 \inf_{t\in [0, T_*]}\inf_{x\in \Rr^d}   \mathrm{RT}(\eta(t)) \ge \tdm\ma.
 \eq
 \end{lemm}
 \begin{proof}
 Set
 \[
 K=\lbb\mathfrak B(\eta_n)J(\eta_n)\rbb\eta_n-\lbb\mathfrak B(\eta)J(\eta)\rbb\eta.
 \]
 Arguing as in the proof of \eqref{A(t):Hs-tdm} we find that
\[
 \begin{aligned}
\| K\|_{H^{s-2}}& \leq \mathcal F(\|(\eta_n, \eta) \|_{H^s})\|\eta_n-\eta\|_{H^{s-1}}\\
& \les \mathcal F(\|\eta_0\|_{H^s}, \ma^{-1})\|\eta_n-\eta\|_{H^{s-1}}.
\end{aligned}
\]
On the other hand, by estimating each term in $K$  we have
\[
\| K\|_{H^{s-1}}\leq \mathcal F(\|(\eta_n, \eta) \|_{H^s}) \les \mathcal F(\|\eta_0\|_{H^s}, \ma^{-1}).
\]
Choosing $s'\in (\max\{\frac{d}{2}, s-2\}, s-1)$, then interpolating the above estimates gives 
\[
\| K\|_{L^\infty([0, T_*]; L^\infty(\Rr^d)}\le \mathcal F(\|\eta_0\|_{H^s}, \ma^{-1})\|\eta_n-\eta\|_{L^\infty([0, T]; H^{s-1})}^\tt
\]
for some $\tt\in (0, 1)$. Then, \eqref{bound:limit3} follows from this and \eqref{RT:etas}.

 \end{proof}
 Now in view of the properties \eqref{bound:limit1}, \eqref{bound:limit2} and \eqref{bound:limit3} of $\eta$, we see that in the proof of  \eqref{proof:finalest}, if we replace $\eta^{(\s_1)}$ with $\eta_n$, $\eta^{(\s_2)}$ with $\eta$, and $(\s_1, \s_2)$ with $(\s_n, 0)$, then we obtain the convergence estimate 
\bq\label{proof:conv1}
\|\eta_n-\eta\|_{L^\infty([0,T_*];H^{s-1})}+\|\eta_n-\eta\|_{L^2([0,T_*];H^{s-\mez})} \le  \sqrt{\s_n}\mathcal F(\|\eta_0\|_{H^s}, \ma^{-1}).
\eq
Furthermore, assume that $\eta_1$ and $\eta_2$ are two solutions on $[0, T_*]$ of the Muskat problem without surface tension with the same initial data $\eta_0$ and that both satisfy \eqref{bound:limit1}, \eqref{bound:limit2} and \eqref{bound:limit3}. Then the proof of \eqref{proof:finalest} with $\s_1=\s_2=0$ yields that $\eta_1\equiv \eta_2$ on $[0, T_*]$. This proves the uniqueness of $\eta$. In other words, we have obtained an alternative proof for  the local well-posedness of the Muskat problem without surface tension for any subcritical data satisfying \eqref{RT:unilw} and \eqref{sep:initial}. 

The next proposition improves the rate in \eqref{proof:conv1} to the optimal rate.
\begin{prop}\label{prop:conv2}
If in addition  $s\ge 2$, then 
\bq\label{conv:opt}
\begin{aligned}
&\|\eta_n-\eta\|_{L^\infty([0,T_*];H^{s-2})}+ \|\eta_n-\eta\|_{L^2([0,T_*];H^{s-\tdm})} \le   \s\mathcal F(\| \eta_0\|_{H^s}, \ma^{-1}),
\end{aligned}
\eq
where  $\cF:\Rr^+\times \Rr^+\to \Rr^+$ depends only on $(h, s, \mu^\pm, \g)$.
\end{prop}
\begin{proof}
We follow the notation in the proof of Proposition \ref{prop:conv} but set $\eta_1=\eta_n$ and $\eta_2=\eta$. Then, $\eta_\delta=\eta_n-\eta$ satisfies \begin{equation} \label{deltSys2}
(\mu^+ + \mu^-) \partial_t \eta_\delta = -\mathfrak g (L_\delta \eta_1+ L_2\eta_\delta)- \mathfrak s_n L_1 H(\eta_1).
\end{equation}
 Applying  \eqref{continuity:JGL} and \eqref{boundH}  with $\sigma=s-\frac32\ge \mez$ yields
\bq\label{L1H1:2}
\| L_1 H(\eta_1)\|_{H^{s-\frac52}}\le \cF(\|\eta_1\|_{H^s})\| H(\eta_1)\|_{H^{s-\tdm}}\le  \cF(\|\eta_1\|_{H^s})\| \eta_1\|_{H^{s+\mez}}.
\eq
Next we paralinearize $L_2$ and $L_\delta$. For $L_2$ we apply  \eqref{Lprincipal} with $\sigma=s-\tdm\ge \mez$
\bq\label{lin:L2:2}
L_2 \eta_\delta = T_{\lambda_2}\eta_\delta+  O_{H^{s-\tdm} \to H^{s -\frac52+\delta}}(\mathcal F(N_s))\eta_\delta.
\eq
$L_\delta$ can be written as in \eqref{expand:Ldelta}. Using \eqref{DNprincipal} with $\sigma = s-\tdm\ge \mez$ together with the fact that $J _\delta^--J^+_\delta=0$, we obtain
\begin{align*}
&\sum_{\pm = +,-} G_2^\pm  J _\delta^\pm =  \sum_{\pm = +,-} O_{\wt H^{s-\tdm}_\pm\to H^{s-\frac52+\delta}}(\mathcal F(N_s)) J _\delta^\pm.
\end{align*}
Applying  Proposition \ref{gamContLem} with $\sigma=s-\tdm\ge \mez$, we obain
 \begin{align*}
 J^\pm_\delta=O_{H^s\to \wt H^{s-\tdm}_\pm}(\| \eta_\delta\|_{H^{s-\tdm}}),
 \end{align*}
and hence
\begin{align*}
\sum_{\pm = +,-} G_2^\pm  J _\delta^\pm =   O_{ H^s\to H^{s-\frac52+\delta}}(\mathcal F(N_s)\|\eta_\delta\|_{H^{s-\tdm}}).
\end{align*}
On the other hand,  Theorem \ref{DNcont}  can be applied with $\sigma=s-\tdm\ge \mez$, implying 
\[
\sum_{\pm = +,-} G_\delta^\pm  J _1^\pm f=-T_{ \lambda_1 \lbb \mathfrak B_1  J _1\rbb f}\eta_\delta  -T_{\lbb \mathfrak V_1  J _1\rbb f} \cdot\nabla  \eta_\delta + O_{H^s \to H^{s-\frac52+\delta}}\big(\mathcal F(N_s) \|\eta_\delta\|_{H^{s-\tdm}}\big)f.
\]
We thus obtain 
\bq\label{Ldelt:2}
\begin{aligned}
L_\delta f &=- T_{ \lambda_1 \lbb \mathfrak B_1  J _1\rbb f}\eta_\delta  -T_{\lbb \mathfrak V_1  J _1\rbb f} \cdot\nabla  \eta_\delta  +O_{ H^s\to H^{s-\frac52+\delta}}(\mathcal F(N_s)\|\eta_\delta\|_{H^{s-\tdm}})f.
\end{aligned}
\eq
Applying this with $f=\eta_1$, then combining with \eqref{lin:L2:2} and symmetrizing we arrive at
\begin{align*}
L_2\eta_\delta+L_\delta \eta_1&= T_{ (\lambda (1-\lbb \mathfrak B  J \rbb \eta))_\alpha}\eta_\delta  -T_{(\lbb \mathfrak V  J\rbb \eta)_\alpha} \cdot\nabla  \eta_\delta \\
&\qquad + O_{H^{s-\tdm} \to H^{s -\frac52+\delta}}(\mathcal F(N_s))\eta_\delta+O_{ H^s\to H^{s-\frac52+\delta}}(\mathcal F(N_s)\|\eta_\delta\|_{H^{s-\tdm}})\eta_1.
\end{align*}
Plugging this and \eqref{L1H1:2} into \eqref{deltSys} leads to
\bq\label{proof:lin1:2}
\begin{aligned}
&(\mu^+ + \mu^-) \partial_t \eta_\delta =- \g T_{ (\lambda (1-\lbb \mathfrak B  J \rbb \eta))_\alpha}\eta_\delta  +\g T_{(\lbb \mathfrak V  J\rbb \eta)_\alpha} \cdot\nabla  \eta_\delta+\cR'_1+\cR_2',\\
&\| \cR_1'\|_{H^{s-\frac52+\delta}}\le \g\cF(N_s)\|\eta_\delta\|_{H^{s-\tdm}},\\
&\| \cR_2'\|_{H^{s-\frac52}}\le\s_n\cF(N_s) \| \eta_1\|_{H^{s+\mez}}.
\end{aligned}
\eq
Next we set  $\eta_{\delta, s-2}=\langle D\rangle^{s-2}\eta_\delta$ and commute the first equation in \eqref{proof:lin1:2} with $\langle D\rangle^{s-2}$ to obtain after applying Theorem \ref{theo:sc} that
\bq\label{proof:lin3:2}
\begin{aligned}
&(\mu^+ + \mu^-) \partial_t \eta_{\delta, s-2} =- \g T_{ (\lambda (1-\lbb \mathfrak B  J \rbb \eta))_\alpha}\eta_{\delta, s-2}  +\g i\mathrm{Re}\big(T_{(\lbb \mathfrak V  J\rbb \eta)_\alpha \cdot\xi}\big)  \eta_{\delta, s-2}+\cR_1+\cR_2,\\
&\| \cR_1\|_{H^{-\mez+\delta}}\le \g\cF(N_s)\|\eta_{\delta, s-2}\|_{H^{\mez}},\\
&\| \cR_2\|_{H^{-\mez}}\le\s_n\cF(N_s) \| \eta_1\|_{H^{s+\mez}},
\end{aligned}
\eq
where $\cF$ depends only on $(h, s, \mu^\pm)$. An $L^2$ energy estimate as in \eqref{eest:etadelta} yields
\bq\label{eest:etadelta:2}
\begin{aligned}
&\frac{(\mu^+ + \mu^-)}{2}\frac{d}{dt}\|\eta_{\delta, s-2} \|_{L^2}^2  +\frac{\g}{\mathcal{Q}_{T_*}}\| \eta_{\delta, s-2} \|^2_{H^\mez}\\
&\quad\le \g\mathcal{Q}_{T_*}\|\eta_{\delta, s-2}\|_{H^{\mez-\delta}}\| \eta_{\delta, s-2}\|_{H^{\mez}}+\s_n\mathcal{Q}_{T_*} \| \eta_1\|_{H^{s+\mez}}\| \eta_{\delta, s-2}\|_{H^{\mez}},
\end{aligned}
\eq
where $\mathcal{Q}_{T_*}$ is given by \eqref{Q:conv}. Interpolating as in \eqref{interpolate:conv} we obtain 
\bq\label{diff:ineq1:2}
\begin{aligned}
&\frac{(\mu^+ + \mu^-)}{2}\frac{d}{dt}\|\eta_{\delta, s-2} \|_{L^2}^2  +\frac{\g}{\mathcal{Q}_{T_*}}\| \eta_{\delta, s-2}\|_{H^{\mez}}^2\le \g\mathcal{Q}_{T_1}\| \eta_{\delta, s-2}\|^2_{L^2}+\s_n^2\mathcal{Q}_{T_*} \| \eta_1\|_{H^{s+\mez}}^2.
\end{aligned}
\eq
From the uniform estimate \eqref{regularity:etas} we have
\[
\int_0^{T_*}\|\eta_1\|^2_{H^{s+\mez}}dt  \le  \mathcal F(\|\eta_0\|_{H^s}, \ma^{-1}).
\]
Thus, applying G\"onwall's lemma to \eqref{diff:ineq1:2} we arrive at \eqref{conv:opt}.
\end{proof}
\appendix

\section{Paradifferential Calculus}\label{appendix}
In this appendix, we recall the symbolic calculus of Bony's paradifferential calculus. See \cite{Bony, MePise}.
\begin{defi}\label{defi:para}
1. (Paradifferential symbols) Given~$\rho\in [0, \infty)$ and~$m\in\Rr$,~$\Gamma_{\rho}^{m}(\Rr^d)$ denotes the space of
locally bounded functions~$a(x,\xi)$
on~$\Rr^d\times(\Rr^d\setminus 0)$,
which are~$C^\infty$ with respect to~$\xi$ for~$\xi\neq 0$ and
such that, for all~$\alpha\in\Nn^d$ and all~$\xi\neq 0$, the function
$x\mapsto \partial_\xi^\alpha a(x,\xi)$ belongs to~$W^{\rho,\infty}(\Rr^d)$ and there exists a constant
$C_\alpha$ such that,
\begin{equation*}
\forall |\xi|\ge \mez,\quad 
\Vert \partial_\xi^\alpha a(\cdot,\xi)\Vert_{W^{\rho,\infty}(\Rr^d)}\le C_\alpha
(1+|\xi|)^{m-|\alpha|}.
\end{equation*}
Let $a\in \Gamma_{\rho}^{m}(\Rr^d)$, we define the semi-norm
\begin{equation}\label{defi:norms}
M_{\rho}^{m}(a)= 
\sup_{|\alpha|\le 2(d+2) +\rho ~}\sup_{|\xi| \ge \mez~}
\Vert (1+|\xi|)^{|\alpha|-m}\partial_\xi^\alpha a(\cdot,\xi)\Vert_{W^{\rho,\infty}(\Rr^d)}.
\end{equation}
2. (Paradifferential operators) Given a symbol~$a$, we define
the paradifferential operator~$T_a$ by
\begin{equation}\label{eq.para}
\widehat{T_a u}(\xi)=(2\pi)^{-d}\int \chi(\xi-\eta,\eta)\widehat{a}(\xi-\eta,\eta)\Psi(\eta)\widehat{u}(\eta)
\, d\eta,
\end{equation}
where
$\widehat{a}(\theta,\xi)=\int e^{-ix\cdot\theta}a(x,\xi)\, dx$
is the Fourier transform of~$a$ with respect to the first variable; 
$\chi$ and~$\Psi$ are two fixed~$C^\infty$ functions such that:
\begin{equation}\label{cond.psi}
\Psi(\eta)=0\quad \text{for } |\eta|\le \frac{1}{5},\qquad
\Psi(\eta)=1\quad \text{for }|\eta|\geq \frac{1}{4},
\end{equation}
and~$\chi(\theta,\eta)$ 
satisfies, for~$0<\eps_1<\eps_2$ small enough,
$$
\chi(\theta,\eta)=1 \quad \text{if}\quad |\theta|\le \eps_1| \eta|,\qquad
\chi(\theta,\eta)=0 \quad \text{if}\quad |\theta|\geq \eps_2|\eta|,
$$
and such that
$$
\forall (\theta,\eta), \qquad | \partial_\theta^\alpha \partial_\eta^\beta \chi(\theta,\eta)|\le 
C_{\alpha,\beta}(1+| \eta|)^{-|\alpha|-|\beta|}.
$$\end{defi}
\begin{theo}\label{theo:Op}
For all $m\in \Rr$, if $a\in \Gamma^m_0$ then 
\bq
T_a=O_{Op^m}\big(M^m_0(a)\big).
\eq
\end{theo}
\begin{theo}[Symbolic calculus] \label{theo:sc}
Let $a \in \Gamma_r^{m}, a'\in \Gamma_r^{m'}$ and set  $\delta = \min\{1,r\}$.   Then, \\
\noindent (i)  
\begin{align}
T_{a}T_{a'} = T_{aa'} +  O_{Op^{m+m'-\delta}}\Big( M_r^{m}(a)M_0^{m'}(a')+M_0^{m_1}(a)M_r^{m'}(a')\Big);
\end{align}
\noindent (ii) 
\begin{equation}
T_{a}^* = T_{\overline a}+ O_{Op^{m-\delta}}\big(M_r^m(a)\big). 
\end{equation}  
\end{theo}
\begin{rema}\label{rema:low}
In the definition \eqref{eq.para} of paradifferential operators, the cut-off $\Psi$ removes the low frequency part of $u$. In particular, if $a\in \Gamma^m_0$  then 
\[
\Vert T_a u\Vert_{H^\sigma}\le CM_0^m(a)\Vert \nabla u\Vert_{H^{\sigma+m-1}}=CM_0^m(a)\Vert  u\Vert_{H^{\sigma+m, 1}},
\]
and similarly for other estimates involving paradifferential operators.
\end{rema}
\begin{prop}[G\r arding's inequality]
Assume $a\in\Gamma^m_r$ with $m\in \Rr$ and $r\in (0, 1]$ such that for some $c>0$
\bq\label{elliptic:op}
\inf_{(x,\xi)\in \Rr^d\times (\Rr^d\setminus\{0\})} \mathrm{Re}(a(x,\xi)) \ge c|\xi|^m.
\eq
  Then, for all $\sigma \in \mathbb R$, there exists $\mathcal F:\Rr^+\times \Rr^+\to \Rr^+$ nondecreasing such that
\begin{equation}
\|\Psi(D)u\|_{H^\frac{m}{2}}^2\leq \mathcal F(M_r^m(a),c^{-1})\Big(\mathrm{Re}(T_a u,u)_{L^2} + \| u\|_{H^{1, \frac{m-r}{2}}}^2\Big)\label{wgi:0}
\end{equation}
and 
\begin{equation}
\|\Psi(D)u\|_{H^\frac{m}{2}}^2\leq \mathcal F(M_r^m(a),c^{-1})\Big(\mathrm{Re}(T_a u,u)_{L^2} + \|u\|_{H^{\frac{m}{2}}}\| u\|_{H^{1, \frac{m}{2}-r}}\Big)\label{wgi}
\end{equation}
provided that both sides are finite. Here,  $\Psi(D)$ is the Fourier multiplier with symbol $\Psi$  given by \eqref{cond.psi}.
\end{prop}
\begin{proof}
We have 
\begin{align*}
\mathrm{Re}(T_a u,u)_{L^2}&=\mez \big((T_a u,u)_{L^2}+(T^*_a u,u)_{L^2}\Big)\\
&=(T_{\mathrm{Re}(a)} u,u)_{L^2}+\mez((T_a^*-T_{\overline a})u,u)_{L^2}).
\end{align*}
According to Theorem \ref{theo:sc} (ii), $T_a^*-T_{\overline a}$ is of order $m-r$ and 
\[
\| ((T_a^*-T_{\overline a})u,u)_{L^2})\|\le \| (T_a^*-T_{\overline a})u\|_{H^{\frac{m-r}{2}}}\| u\|_{H^{\frac{m-r}{2}}}\le CM^m_r(a) \| u\|^2_{H^{\frac{m-r}{2}}}.
\]
Set $b=(\mathrm{Re}(a))^\mez$. By virtue of \eqref{elliptic:op} we have $b\in \Gamma^{\frac{m}{2}}_r$ and $M^{\frac{m}{2}}_r(b)\le \cF(M^m_r(a))$. We write
\begin{align*}
(T_{\mathrm{Re}(a)} u,u)_{L^2}&=(T_bT_bu, u)_{L^2}+((T_{b^2}-T_bT_b)u, u)_{L^2}\\
&=(T_bu, T_b^*u)_{L^2}+((T_{b^2}-T_bT_b)u, u)_{L^2}\\
&=\|T_bu\|_{L^2}^2+(T_bu, (T_b^*-T_b)u)_{L^2}+((T_{b^2}-T_bT_b)u, u)_{L^2}.
\end{align*}
Applying Theorem \ref{theo:sc} (ii) once again we deduce that $T_b^*-T_b$ is of order $\frac{m}{2}-r$ and 
\[
 \big|(T_bu, (T_b^*-T_b)u)_{L^2}\big|\le \| T_bu\|_{H^{-\frac{r}{2}}} \|(T_b^*-T_b)u)_{L^2}\|_{H^\frac{r}{2}}\le \cF(M^m_r(a), c^{-1}) \| u\|^2_{H^{1, \frac{m-r}{2}}},
\]
where we used Remark \ref{rema:low} in the last inequality. On the other hand, an application of Theorem \ref{theo:sc} (i) yields 
\[
\big|((T_{b^2}-T_bT_b)u, u)_{L^2}\big|\le \cF(M^m_r(a), c^{-1}) \| u\|^2_{H^{1, \frac{m-r}{2}}}.
\] 
Thus,  we obtain
\bq\label{Garding1}
 \|T_bu\|_{L^2}^2\le (T_{\mathrm{Re}(a)} u,u)_{L^2}+\cF(M^m_r(a), c^{-1})\| u\|^2_{H^{1, \frac{m-r}{2}}}.
\eq
By shifting derivative differently in the above inner products, we have the variant
\bq\label{Garding2}
 \|T_bu\|_{L^2}^2\le (T_{\mathrm{Re}(a)} u,u)_{L^2}+\cF(M^m_r(a), c^{-1}) \| u\|_{H^{1, \frac{m}{2}}}\| u\|_{H^{1, \frac{m}{2}-r}}.
\eq
Next we note that $T_{b^{-1}}T_b-\Psi(D)=T_{b^{-1}}T_b-T_1$ is of order $-r$ and 
\bq\label{Garding3}
\begin{aligned}
\| \Psi(D)u\|_{H^\frac{m}{2}}&=\|T_{b^{-1}}T_bu\|_{H^\frac{m}{2}}+\|(T_{b^{-1}}T_b-\Psi(D))u\|_{H^\frac{m}{2}}\\
&\le \cF(M^m_r(a), c^{-1})\| T_bu\|_{L^2}+ \cF(M^m_r(a), c^{-1})\| u\|_{H^{1, \frac{m}{2}-r}}.
\end{aligned}
\eq
Finally, a combination of \eqref{Garding1} and \eqref{Garding3} leads to \eqref{wgi:0}, and a combination of \eqref{Garding2} and \eqref{Garding3} leads to \eqref{wgi}.
\end{proof}
The proof of \eqref{Garding3} also proves the following lemma.
\begin{lemm}\label{lemm:invertOp}
Let $a\in \Gamma^m_r$, $r\in (0, 1]$, be a real symbol satisfying $a(x, \xi)\ge c|\xi|^m$ for all $(x, \xi)\in \Rr^d\times \Rr^d$. Then for all $s\in \Rr$ we have
\bq
\| \Psi(D)u\|_{H^{s}}\le\cF(M^m_r(a), c^{-1})\| T_au\|_{H^{s-m}}+ \cF(M^m_r(a), c^{-1})\| u\|_{H^{1, s-r}}.
\eq
\end{lemm}

\vspace{.1in}
\noindent{\bf{Acknowledgment.}} 
The work of HQN was partially supported by NSF grant DMS-1907776. The authors thank B. Pausader for discussions  about the Muskat problem.  We would like to thank the reviewer for his/her careful reading and helpful suggestions.

\end{document}